\documentclass[11pt]{article}
\usepackage{amssymb}
\usepackage{amsmath}

\newcommand{\asp}{A_{N-(p-1)}^{p-1}}
\newcommand{\av}{\mbox{{\bf Av}}}
\newcommand{\be}{\mbox{{\bf E}}}

\newcommand{\entk}{{\cal E}_{n,k,t}}

\newcommand{\eti}{\eta_{J}}

\newcommand{\gip}{g_{i_1,\ldots,i_p}}

\newcommand{\ou}{[0,1]}

\newcommand{\sii}{\sigma_i}
\newcommand{\siu}{\sigma_{i_1}}
\newcommand{\sip}{\sigma_{i_p}}

\newcommand{\sbn}{\sum_{J\in Q_{N,k}^p}}
\newcommand{\sln}{\sum_{l=1}^n}
\newcommand{\slp}{\sum_{1\le l<l'\le n}}

\newcommand{\ssm}{\Sigma_{N-1}}
\newcommand{\ssn}{\Sigma_N}
\newcommand{\J}{\hat J}

\newcommand{\R}{\mathbb R}
\newcommand{\N}{\mathbb N}

\newcommand{\ce}{\mathcal E}

\newcommand{\al}{\alpha}
\newcommand{\ep}{\varepsilon}

\newcommand{\si}{\sigma}

\newcommand{\vp}{\varphi}

\newcommand{\lp}{\left(}
\newcommand{\rp}{\right)}
\newcommand{\lc}{\left[}
\newcommand{\rc}{\right]}
\newcommand{\lcl}{\left\{}
\newcommand{\rcl}{\right\}}
\newcommand{\lla}{\left\langle}

\newcommand{\rak}{\right\rangle_{k}}

\newtheorem{theorem}{Theorem}[section]

\newtheorem{corollary}[theorem]{Corollary}

\newtheorem{definition}[theorem]{Definition}

\newtheorem{lemma}[theorem]{Lemma}

\newtheorem{proposition}[theorem]{Proposition}
\newtheorem{remark}[theorem]{Remark}

\begin{document}
 \thispagestyle{empty}
 \begin{center}
 {\Large\bf The p-spin interaction model}\\
{\Large\bf with external field}

\vspace{0.3cm}

 by

\vspace{0.3cm}

{\bf Xavier Bardina}\footnote{ Partially supported
by DGES grants BFM2000-0009, BFM2000-0607,  HF2000-0002.}, 
{\bf David M\'arquez-Carreras$^2$, 
Carles Rovira}\footnote{{Partially supported by DGES grants
BFM2000-0607 and HF2000-0002
.\hfill}}\\  
and\\ 

{\bf Samy
Tindel$^3$}

\vspace{0.1cm}

$^1$
{\it Departament de Matem\`atiques, Universitat Aut\`onoma de Barcelona}
\\
\it 08193 Bellaterra, Barcelona, Spain
\\
{\it e-mail: bardina@mat.uab.es}
\\
 $^2$
{\it Facultat de Matem\`atiques,
Universitat de Barcelona,}
\\
\it Gran Via 585,
08007-Barcelona,
Spain
\\
{\it e-mail: marquez@mat.ub.es, rovira@mat.ub.es}
\\
$^3$ {\it D\'epartement de Math\'ematiques,
Institut Galil\'ee - Universit\'e Paris 13,}
\\
{\it Avenue  J. B. Cl\'ement,
93430-Villetaneuse,
France}
\\
{\it e-mail: tindel@math.univ-paris13.fr}

 \end{center}

\vspace{0.3cm}

\begin{abstract}
This paper is devoted to a detailed 
study of a $p$-spins interaction model 
with external field, including some sharp 
bounds on the speed of self averaging 
of the overlap as well as a central limit 
theorem for its fluctuations, the thermodynamical 
limit for the free energy and the definition of 
an Almeida-Thouless type line. Those results show 
that the external field dominates the tendency 
to disorder induced by the increasing level 
of interaction between spins, and 
our system will share many of its features
with the SK model, which is certainly not the case when 
the external magnetic field vanishes.
\end{abstract}

\newpage

\section{Introduction}

The high temperature regime of the 
Sherrington-Kirkpatrick model of spin glasses, 
with or without external field, is now understood 
in many of its essential features: the overlap $R_{1,2}$ 
of two configurations has been shown to be a 
central object of study for the whole system 
(see \cite{MPV}), the thermodynamical limit of 
$R_{1,2}$ and of the free energy $Z_N$ have been 
computed (see e.g. \cite{T98}), and a number of Central 
Limit Theorems for the fluctuations of those 
quantities have also been established in different 
contexts (see \cite{ALR}, \cite{CN}, \cite{Tbk}), 
giving a rather complete picture of the model.

On the other hand, the results concerning a natural 
generalization of the SK model, namely the $p$-spins 
interaction model, are scarce (see however \cite{Tp} 
on the low temperature regime and \cite{BKL} for some 
fluctuation results for the free energy), especially 
when an external field is considered. The purpose of 
the present paper is then to fill this gap: we will 
consider a spin glass model, whose configuration space 
is $\ssn=\{ -1,1 \}^N$. Let $\mu_N$ be the uniform 
measure on $\ssn$. The energy of a given configuration 
$\si\in\ssn$ will be represented by a Hamiltonian $H(\si)$, 
and we are concerned with the Gibbs measure $G=G_N$, 
whose density with respect to $\mu_N$ is $Z_N^{-1}e^{- H}$, 
where $Z_N$ is the normalization factor
$$
Z_N=\sum_{\si\in\ssn}\exp\lp- H(\si)  \rp.
$$ 
The Hamiltonian under consideration here will be defined by
$$
-H_{N,\beta,h}(\si)=\beta u_N \sum_{(i_1,\dots,i_p) 
\in A_N^p}\gip \siu\ldots\sip
+h\sum_{i=1}^{N}\sii,
$$
with
\begin{eqnarray}
u_N&=&\lp \frac{p!}{2N^{p-1}} \rp^{\frac12},\nonumber\\
A_N^p&=&\lcl (i_1,\ldots,i_p)\in\mathbb{N}^p; 1\le i_1<\cdots<i_p\le N  
\rcl,\nonumber
%\label{edefa}
\end{eqnarray}
where  the parameter $\beta$ represents the inverse of the temperature
and where $g=\{ \gip;(i_1,\ldots,i_p)\in A_N^p \}$ is a family 
of independent standard Gaussian random variables.
The strictly positive parameter $h$ stands for the external magnetic 
field, 
under which the spins tend to take the same value $+1$.
We will denote by $\langle f\rangle$ the average of a function 
$f:\Sigma_N\to\R$ with respect to $G_N$, as well as the average 
of a function $f:\Sigma_N^n\to\R$ with respect to $G_N^{\otimes n}$, 
without mentionning the number $n$ of independent copies of the spins 
configurations, i.e.
$$
\langle f\rangle=Z_N^{-n}\sum_{(\si^1,\ldots,\si^n)\in\Sigma_N^n}
f(\si^1,\ldots,\si^n)\exp\lp -\sum_{l\le n}H_{N,\beta,h}(\si^l) \rp.
$$
We write $\nu(f)=\be \langle f\rangle$.
Our aim here is then to give a detailed account 
on the limiting behavior of this system when 
$N\to\infty$, when $\beta$ is bounded from above by a constant 
$\beta_p$.

Notice that some of the features of 
the SK model are shared by our $p$-spins 
interaction model. For instance, the study 
of the overlap of two configurations, defined by
$$
R_{1,2}=\frac{1}{N}\sum_{i\le N}\sii^1\sii^2,
$$ 
where $\si^1,\si^2$ are understood as two independent 
configurations under $G_N$, will be again one of the 
main steps to understand the limiting behavior of the 
system, though it generally appears under the form 
$R_{1,2}^{p-1}$ (for instance in 
our first occurrence of the cavity method, yielding 
Proposition \ref{pderivnu})
, leading to some technical complications. 
Our first result will then be to show that, 
for $\beta$ small enough, $R_{1,2}$ will self 
average into a constant $q=q_p$, implicitely given 
as the unique solution to
$$
q=\be \Bigg[ \tanh^2\Bigg(\beta \Big(\frac{p}{2}\Big)^{\frac{1}{2}}
q^{\frac{p-1}{2}} Y +h\Bigg)\Bigg],
$$
where $Y$ stands for a standard Gaussian random variable. 
In particular, it will be easily shown that $q_p$ will 
tend to $\tanh^2(h)$ as $p$ grows to 
$\infty$, showing that the natural tendency to the disorder 
induced by the increasing level of interaction between the 
spins will be dominated by the presence of the external field $h$.

It will be natural then to obtain some extra information on the 
exponential moments of $N(R_{1,2}-q)$, from which we will 
be able to get the estimate
\begin{equation}\label{bndrmq}
\nu\lp R_{1,2}-q \rp\le\frac{L(p,\beta)}{N},
\end{equation}
giving a sharp bound on the speed of self averaging of 
$R_{1,2}$. All those considerations on the overlap will 
yield the following replica-type formula for $Z_N$:
\begin{multline*}
\lim_{N\uparrow \infty} p_N(\beta,h,p)=
\frac{\beta^2}{4}\big[1-pq^{p-1}+(p-1)q^p\big]\\
+ \log 2 + \be\Bigg[\log \cosh
\Big[\beta \Big(\frac{p}{2}\Big)^{\frac{1}{2}}
q^{\frac{p-1}{2}}Y+h\Big]\Bigg].
\end{multline*}
Some further computations on the second moments 
of $R_{1,2}^{p-1}$  will then lead us to the 
definition of an Almeida-Thouless line, which 
should give the limit of the high temperature 
region for our model, and is defined by
$$
1-\frac{\beta^2p(p-1)q^{p-2}}{2}\be\lc \cosh^{-4}
\lp\Big(\frac{p}{2}\Big)^{\frac{1}{2}}
q^{\frac{p-1}{2}} Y +h\rp \rc >0.
$$
Notice that the fact that $q_p\to\tanh^2(h)$ when $p\to\infty$
will immediately imply that, if we denote by $\beta_{\mbox{\tiny at}}$
the boundary of this Almeida-Thouless line, then 
$\beta_{\mbox{\tiny at}}\to\infty$ when $p\to\infty$ (see Remark \ref{limbeat}).

Our last result will be a central limit theorem for $R_{1,2}$: 
we will show that, for the  typical disorder $g$, 
the quantity $N^{1/2}(R_{1,2}^{p-1}-q^{p-1})$ will 
converge to a Gaussian random variable whose variance 
will be identified explicitely. Notice that this behavior 
is quite different from the picture given by \cite{BKL}. 
Indeed, when $h=0$, the rate of fluctuation of $Z_N$ is 
shown to be of order $N^{(p-2)/2}$, increasing thus with the 
number of interactions. In our case, the presence of 
the external field $h$ will stabilize the behavior 
of the self averaging, which will occur at the same 
speed as in the SK case.

Of course, our methods of proofs are much indebted 
to the great influence of \cite{Tbk}, through 
the rigorous introduction of the cavity method 
as well as for some key ideas for further 
computations of moments and limit theorems. 
However, the presence of an increasing number 
of interactions requires a careful analysis 
of the different quantities considered at 
each step of our calculations, especially 
in the identification of all the negligible 
terms involved. This is why we include almost 
all the details of the computations in our proofs, 
which, we hope, will make the lecture of the paper 
easier, though certainly cumbersome.

Our paper is organized as follows: 
at section 2, we will give some preliminary 
results on the cavity method for the $p$-spin model, 
allowing to reduce our system of size $N$ into a system 
of size  $N-k$ for arbitrary $1 \le k \le N$. Section 3 is 
devoted to a preliminary study of $R_{1,2}$, 
including the self averaging result, the existence 
of exponential moments, and the bound (\ref{bndrmq}). 
Section 4 will then give the limiting behavior of 
$\frac1N\be[\log(Z_N)]$. Section 5 will focus on 
the definition of the Almeida-Thouless type line, 
while at Sections 6 and 7 we will establish the CLT for $R_{1,2}$. 
Finally, in the Appendix we recall the definitions of all
the sets appearing througout the paper.

In the sequel, the size of a given finite set $D$ will be denoted by $|D|$. 
Troughout the paper, $P_m(N)$ denotes a  polynomial of order $m$  in $N$.
We will also denote by K almost all the constants, although
their value may change from line to line. We will omit their dependence on
$k$ (the size of the cavity we will create) and $n$ (the number of copies of $G_N$ considered).

\section{The cavity method}
In this Section, we will introduce one of the basic tools we will use all along the paper, namely the $k$-cavity method, that allows to quantify in a certain way the difference between our original system and a system where the $k$ last spins are independent from the other ones. We will first introduce the basic notations we will need further on, then get some general results for the $k$-cavity, and eventually a simplified version of some of these results for the particular case of a 1-cavity.

\subsection{Notations and definitions}

For $k\in \{1,\dots,N-1\}$ and $\beta>0$, let  
$$
%\begin{equation}\label{defbetaprime}
\beta_k=\lp \frac{N-k}{N} \rp^{\frac{p-1}{2}}\beta,
$$
%\end{equation}
that will play the role of $\beta$ for our reduced system.
Define  the following set:
$$Q_{N,k}^p= \big\{J=(i_1,\dots,i_p)\in \mathbb{N}^p;1 \le  i_1<\cdots <i_p \le N, i_p>N-k\big\},$$
and, for $J\in Q_{N,k}^p$, set $m=\max\{j, i_j\le N-k\}$, and let $I,I^c$ be defined by
\begin{equation}\label{estar}
I=(i_1,\dots,i_m),\qquad 
I^c=(i_{m+1},\dots,i_p).
\end{equation}
Observe that we should write  $I=I(J), I^c=I^c(J)$, but we will omit this dependence for sake of readability.
Using these sets, we define
$$
\eta_J=\prod_{i_j\in I} \si_{i_j},\qquad\quad
\ep_J=\prod_{i_j\in I^c} \si_{i_j},
$$
and
\begin{eqnarray*}
g_{(k)}(\eta,\ep)&=&\beta u_N \sum_{i_p>N-k} g_{i_1,\dots,i_p}\si_{i_1}\cdots\si_{i_p}\\
&=&\beta u_N \sum_{J\in Q_{N,k}^p} g_J\ \eta_J\ \ep_J.
\end{eqnarray*}
The basic idea of the $k$-cavity method is to regroup the Hamiltonian as follows:
$$
-H_{N,\beta,h}(\si)=-H_{N-k,\beta_k,h}(\si)+ g_{(k)}(\eta,\ep)+
h\sum_{i=1}^k \ep_i,
$$
with $\ep_i=\si_{N-i+1}$. We will denote then by $\langle \cdot \rangle_k$ the average
with respect to the Gibbs measure on $\Sigma_{N-k}$ relative to the Hamiltonian
$H_{N-k,\beta_k,h}$.
As usual in the spin glasses theory, the cavity method
will become a powerful tool through a construction of a continuous path between
the original configuration, and a configuration where the $k$ last spins are independent
of the others.
Set then, for $t\in [0,1]$ and a constant $q\in [0,1]$ to be precised later,
\begin{equation}\label{defpath}
g_{(k),t}(\eta,\ep)=t^{\frac{1}{2}}g_{(k)}(\eta,\ep)+\beta u_N q^{\frac{p-1}{2}}(1-t)^{\frac{1}{2}}
\sum_{J\in Q_{N,k}^p} z_J\ \ep_J,
\end{equation}
where $\{z_J ;J\in  Q_{N,k}^p \}$ 
is a family of independent standard Gaussian random variables, 
also independent of all the disorder $g$. 

Let $n\ge 1$ and $\si^1,\dots,\si^n$ be $n$ independent copies of a $N$-spins
configuration. Let us write
\begin{eqnarray*}
\ce_{n,k,t}&=&\exp\Big\{ \sum_{l=1}^{n}\big( g_{(k),t}(\eta^l,\ep^l)+ h\sum_{i=1}^k \ep_i^l\big)\Big\},\\
Z_{(k),t}&=&\lla\av \ce_{1,k,t}\rak,
\end{eqnarray*}
where $\av$ means average over $\{\ep_i^l=\pm 1, i=1,\dots,k, l=1,\dots,n\}$. Then,
for $f:\ssn^n\longrightarrow \mathbb{R}$, we define
\begin{eqnarray*}
\langle f \rangle_{k,t}&=&\frac{\langle\av f \ce_{n,k,t} \rangle_{k}}{Z_{(k),t} ^n},\\
\nu_{k,t}(f)&=& \be\langle f \rangle_{k,t}.
\end{eqnarray*}
Observe that $\nu(f)=\nu_{k,1}(f)$ for any $k$.

The idea of what are going to state is that  $\nu_{k,0}(f)$ (or a 
slight modification of this quantity) should be simpler to compute than $\nu_{k,1}(f)$ in some interesting cases of functions $f$.
On the other hand, we will relate these two quantities by means of 
\begin{equation}\label{difnuk}
\nu_{k,1}(f)-\nu_{k,0}(f)=\int_0^1 \frac{d}{dt}\nu_{k,t}(f)\ dt,
\end{equation}
or higher-order versions. We will generally write
$$\nu_{k,t}'(f)=\frac{d}{dt}\nu_{k,t}(f).
$$

\subsection{The $k$-cavity}
We will give here some basic relations, allowing to estimate quantities like (\ref{difnuk}) in great generality.
First, we will compute the derivative of $\nu_{k,t}(f)$ with respect to
$t$ in the following way:
\begin{proposition}\label{pderivnu}
For $t\in\ou$ and $f:\ssn^n\to\R$, we have
\begin{eqnarray}\label{derivnu}
\nu_{k,t}'(f)&=&\beta^2u_N^2 \displaystyle\sum_{J\in Q_{N,k}^p}
\Bigg[\nu_{k,t} \Big( f\slp(\eta_{J}^{l}\eta_J^{l'}-q^{p-1})\
\ep_J^{l}\ep_J^{l'}\Big)\nonumber\\
&&-n\nu_{k,t} \Big( f\sum_{l=1}^n(\eta_J^{l}\eta_J^{n+1}-q^{p-1})\
\ep_J^{l}\ep_J^{n+1}\Big)\nonumber\\
&&+\frac{n(n+1)}{2}
\nu_{k,t}\Big( f(\eta_J^{n+1}\eta_J^{n+2}-q^{p-1})\ 
\ep_J^{n+1}\ep_J^{n+2}\Big)\Bigg].\end{eqnarray}
\end{proposition}

%\vspace{0.2cm}

\noindent
{\bf Proof:} 
This proof is an extension of \cite[Proposition 2.4.2]{Tbk}, 
whose details are given for sake of completeness. We have
$$
\frac{d\entk}{dt}=
\Bigg[ \frac{1}{2t^{\frac{1}{2}}}\sln g_{(k)}(\eta^l,\ep^l)
-\frac{\beta u_N q^{\frac{p-1}{2}}}{2(1-t)^{\frac{1}{2}}}
\sbn z_J \sln\ep^l_J \Bigg] \entk,
$$
and hence $\nu_{k,t}'(f)=A_1-A_2$, with
\begin{eqnarray*}
A_1=\frac{1}{2t^{\frac{1}{2}}}&\be\Bigg[&  Z_{(k),t}^{-n}
\sln\lla\av f g_{(k)}(\eta^l,\ep^l)\ \entk \rak\\
&&-n Z_{(k),t}^{-(n+1)}\lla\av f g_{(k)}(\eta^{n+1},\ep^{n+1})\ \ce_{n+1,k,t} \rak
\Bigg],
\end{eqnarray*}
and
\begin{eqnarray*}
A_2=\frac{\beta u_N q^{\frac{p-1}{2}}}{2(1-t)^{\frac{1}{2}}}
\sbn&\be\Bigg[&  Z_{(k),t}^{-n}\sln\lla\av f z_J\ \ep_J^l\ \entk \rak\\
&&-n Z_{(k),t}^{-(n+1)}\lla\av f z_J \ \ep_J^{n+1} \ \ce_{n+1,k,t} \rak
\Bigg].
\end{eqnarray*}
Notice that, in order to obtain the last formula, we have used the basic fact that,  for $\phi:\ssn^{m}\to\R$ and $\hat\phi:\ssn^{\hat m}\to\R$, we have
\begin{multline*}
\lla \phi\lp\si^1,\ldots,\si^{m}  \rp\rak
\lla \hat\phi\lp\si^1,\ldots,\si^{\hat m}  \rp\rak\\
=\lla \phi\lp\si^1,\ldots,\si^{m}  \rp
\hat\phi\lp\si^{m+1},\ldots,\si^{m+\hat m}  \rp\rak.
\end{multline*}
Let us first study $A_2$: integrating 
this expression with respect to $z_J$, and invoking 
the fact that $\be[z F(z)]=\be[F'(z)]$ for a standard 
Gaussian random variable $z$, we get
$$
\partial_{z_J}\entk=
\beta u_N q^{\frac{p-1}{2}} (1-t)^{\frac{1}{2}}
\lp \sum_{l'=1}^n\ep^{l'}_J \rp\ \entk,
$$
and hence
\begin{eqnarray*}
A_2=\frac{\beta^2u_N^2q^{p-1}}{2}\sbn&\be\Bigg[&  
Z_{(k),t}^{-n}\sum_{l=1}^n\sum_{l'=1}^n\lla \av f \ep_J^{l} \ep_J^{l'}\ \entk  \rak\\
&&-n Z_{(k),t}^{-(n+1)}\sln\lla \av f \ep_J^{l} \ep_J^{n+1}\ \ce_{n+1,k,t}  \rak\\
&&-n Z_{(k),t}^{-(n+1)}\sum_{l=1}^{n+1}\lla \av f \ep_J^{l} \ep_J^{n+1}\ \ce_{n+1,k,t}  \rak\\
&&+\kappa_n Z_{(k),t}^{-(n+2)}\lla \av f \ep_J^{n+1} \ep_J^{n+2}\ \ce_{n+2,k,t}  \rak
\Bigg],
\end{eqnarray*}
with $\kappa_n=n(n+1)$. Note that, for any $l\in\{1,\dots,n\}$,
$$
Z_{(k),t}^{-n}\lla\av f \ep_J^{l} \ep^{l}_J\ \entk  \rak
= Z_{(k),t}^{-n} \lla \av f\ \entk  \rak,
$$
and,  since $f$ depends only on $(\si^1,\ldots,\si^n)$, we have
\begin{eqnarray*}
n Z_{(k),t}^{-(n+1)}\lla \av f \ep_J^{n+1} \ep_J^{n+1}\ \ce_{n+1,k,t}  \rak
&=&
n Z_{(k),t}^{-(n+1)}\lla \av f\  \ce_{n+1,k,t}  \rak\\
&=&
n Z_{(k),t}^{-n}\lla \av f \ \ce_{n,k,t}  \rak.
\end{eqnarray*}
Thus
\begin{eqnarray*}
A_2=\beta^2 u_N^2 q^{p-1}\sbn&\be\Bigg[&
Z_{(k),t}^{-n}\sum_{1\le l<l'\le n}
\lla \av f \ep_J^{l} \ep_J^{l'}\ \entk  \rak\\
&&-n Z_{(k),t}^{-(n+1)}\sln\lla \av f \ep_J^{l} \ep_J^{n+1}\ \ce_{n+1,k,t}  \rak\\
&&+\frac{\kappa_n}{2} Z_{(k),t}^{-(n+2)}\lla \av f \ep_J^{n+1} \ep_J^{n+2}\ \ce_{n+2,k,t}  \rak
\Bigg],
\end{eqnarray*}
where we recall that $\kappa_n=n(n+1)$. This can be read as
\begin{eqnarray*}
A_2=\beta^2u_N^2q^{p-1}&\displaystyle\sbn\Bigg[&
\nu_{k,t} \lp f\sum_{1\le l<l'\le n} \ep_J^{l} \ep_J^{l'} \rp 
-n\nu_{k,t}\lp f\sum_{l=1}^n \ep^{l}_J \ep^{n+1}_J  \rp\\
&&+\frac{n(n+1)}{2}\nu_{k,t}\lp f \ep^{n+1}_J \ep^{n+2}_J \rp \Bigg].
\end{eqnarray*}

The same kind of computations can be lead for $A_1$, 
integrating first by parts with respect to the variables $g_I$. In this case, we obtain
\begin{eqnarray*}
A_1=\beta^2u_N^2&\displaystyle\sbn
\Bigg[& \nu_{k,t}\lp f\sum_{1\le l<l'\le n}  \eti^{l} \eti^{l'}\ \ep_J^{l} \ep_J^{l'} \rp \\
&&-n\nu_{k,t}\lp f \sum_{l=1}^n \eti^{l} \eti^{n+1}\ \ep_J^{l} \ep_J^{n+1}  \rp\\
&&+\frac{n(n+1)}{2}\nu_{k,t}\lp f \eti^{n+1} \eti^{n+2}\ \ep_J^{n+1} \ep_J^{n+2} \rp
\Bigg].
\end{eqnarray*}
Substracting $A_1$ and $A_2$, we get the desired result.
\hfill $\Box$

%\vspace{0.2cm}

As a consequence of the last proposition, we can bound $\nu_{k,t}(f)$ 
by $\nu_{k,1}(f)$ as follows:
\begin{proposition}\label{posmaj}
Let $f:\ssn^n\to\R_+$ be a non-negative function. Then, for $N$ large enough, we have
\begin{equation}\label{majnutnu}
\nu_{k,t}(f)\le \exp\big\{ 2\beta^2n^2p(k+1) \big\}\nu(f)
\end{equation}
\end{proposition}

%\vspace{0.2cm}

\noindent
{\bf Proof:}
Appealing to relation (\ref{derivnu}), we obtain, 
for a non-negative $f$, and since $|\eti^l\eti^{l'}-q^{p-1}|\le 2$,
$$
\nu_{k,t}'(f)\ge-4 \beta^2 n^2 u_N^2|Q_{N,k}^p|\nu_{k,t}(f).
$$
Using the expression of $u_N$ and the estimate of Lemma \ref{cardbn} (see the Appendix)
on $|Q_{N,k}^p|$, we get, for a constant $K(p,k)>0$ depending only on $p$
and $k$,
$$
\nu_{k,t}'(f)\ge-2 \beta^2 n^2 p \lp k+\frac{K(p,k)}{N}  \rp \nu_{k,t}(f).
$$
Hence, for $N$ large enough,
$$
\nu_{k,t}'(f)\ge-2 \beta^2 n^2 p (k+1) \nu_{k,t}(f).
$$
Integrating this relation between $t$ and 1, we get
$$
\log\lp \nu(f) \rp-\log\lp \nu_{k,t}(f) \rp
\ge -2 \beta^2 n^2 p (k+1)(1-t)\ge -2 \beta^2 n^2 p (k+1) ,
$$
which yields the announced relation.

\hfill
$\Box$

%\vspace{0.2cm}

We will finish this subsection with a useful result for the $p-1$-cavity. This lemma gives
 an idea of how our computations will become explicit when $\nu_{p-1,0}$ is considered instead of $\nu$.

\begin{lemma}\label{l2.5}
Let $f^-:\Sigma_N^n\to\R$ be a function depending 
on $\{\si_1^l,\ldots,\si_{N-p+1}^{l};$ $l$ $\le n  \}$. 
Let $Y$ be a standard Gaussian random variable.
For $l\le n$, we designate by $\mathcal{M}_l$ an arbitrary subset of $\{1,\ldots,p-1\}$.
Then, for a constant $L_1(\beta)>0$,
\begin{multline*}
\Bigg|
\nu_{p-1,0}\Big(f^-\prod_{l\le n , m_l\in \mathcal{M}_l}\ep_{m_l}^l\Big)-\nu_{p-1,0}(f^-)\\
\times \prod_{j\le p-1}
\be\Bigg[ \tanh^{\sum_{l=1}^n \mathbb{I}_{\{j\in \mathcal{M}_l\}}}
\Big(\beta \lp\frac{p}{2} \rp^{\frac{1}{2}}
q^{\frac{p-1}{2}}Y+h \Big)  \Bigg]
 \Bigg|
\le \frac{L_1(\beta)}{N} \|f^-\|_\infty.
\end{multline*}     
\end{lemma}

\begin{remark} 
In the sequel we use the following notation:
\begin{equation}\label{defqn}
\hat q_n = \be \Big[
\tanh^n\Big(\beta q^{\frac{p-1}{2}} \big(\frac{p}{2}\big)^{\frac{1}{2}}Y
+h\Big)\Big].
\end{equation}
\end{remark}

In order to prove this lemma, it will be useful to 
split $Q_{N,k}^p$ into
\begin{equation}\label{AA}
Q_{N,k}^p= \bar Q_{N,k}^p  \cup \tilde Q_{N,k}^p,
\end{equation}
where
\begin{eqnarray*}
\bar Q_{N,k}^p&=&\big\{J\in  Q_{N,k}^p\ ; m<p-1\big\},\\
\tilde Q_{N,k}^p&=&\big\{J\in  Q_{N,k}^p\ ; m=p-1\big\}.
\end{eqnarray*}
Recall that $m$ is defined in (\ref{estar}) as the maximum of $\{j, i_j\le N-k\}$. Hence, $\bar Q_{N,k}^p$ can also be defined as 
$$\bar Q_{N,k}^p 
= \big\{J=(i_1,\dots,i_p)\in Q_{N,k}^p\ ;\ i_1<\cdots<i_p, i_{p-1}>N-k\big\}.$$

Finally, we need to introduce some additional notation and to give a technical lemma that we will only use in the proof of Lemma \ref{l2.5}, and that expresses the fact that, in (\ref{defpath}), the main part of $\sum_{J \in Q_{N,k}^p} z_J \ep_J$ is given by $\sum_{J \in \tilde Q_{N,k}^p} z_J \ep_J$, which is easier to handle:
let us write
\begin{displaymath}
\hat\ce_{n,k,0}=\exp\Big\{ \sum_{l\le n}
\big( \beta u_N q^{\frac{p-1}{2}}\sum_{i\le k} \hat z_i\ \ep^l_i+ 
h\sum_{i\le k} \ep_i^l\big)\Big\},
\end{displaymath}
where $\{\hat z_i, i=1,\ldots,k\}$ are 
independent zero mean Gaussian random variables with variance
${N-k \choose p-1}$.  For $f\equiv f(\ep_1^l,\ldots,\ep_k^l, l\le n)$, we define
\begin{displaymath}
\hat\nu_{k,0}(f)=\be\Bigg[\frac{\av (f \hat\ce_{n,k,0})}{\hat Z_{(k)}^n}\Bigg],
\end{displaymath}
with $\hat Z_{(k)}=\av \hat\ce_{1,k,0}$.

\begin{lemma}\label{l2.4}
Let $f\equiv f(\ep_1^l,\dots,\ep_k^l, l\le n)$. Then
\begin{displaymath}
|\nu_{k,0}(f)-\hat\nu_{k,0}(f)| \le \frac{K(\beta)}{N} \|f\|_\infty.
\end{displaymath}
\end{lemma}

%\vspace{0.2cm}

\noindent
{\bf Proof:} The arguments are similar to Proposition \ref{pderivnu}. Consider
\begin{displaymath}
\hat\ce_{n,k,t}=\hat\ce_{n,k,0}  \times
\exp\Big\{ \sum_{l=1}^{n}
\beta u_N\ q^{\frac{p-1}{2}}\ t^{\frac{1}{2}} 
\sum_{J \in \bar Q_{N,k}^p} z_{J}\ \ep^l_{J}\Big\},
\end{displaymath}
where
$\{z_{J}, J \in \bar Q_{N,k}^p \}$ are independent standard Gaussian random
variables, and 
define
\begin{displaymath}
\hat\nu_{k,t}(f)=
\be\Bigg[\frac{\av (f \hat\ce_{n,k,t})}{\hat Z_{(k),t}^n}\Bigg],
\end{displaymath}
with $\hat Z_{(k),t}=\av\ \hat\ce_{1,k,t}$.

Note that 
$\hat\nu_{k,1}(f)=\nu_{k,0}(f)$. 
Indeed,
$\hat z_i \sim 
N\big(0,{N-k \choose p-1}\big)$ and thus 
$$\sum_{i=1}^k \hat z_i \ep_i^l
\stackrel{{\rm (d)}}{=}
\sum_{J\in \tilde Q_{N,k}^p} 
z_J \ep_J^l.$$ 
The quantity $\hat\nu_{k,t}(f)$ can be differentiated  once again in $t$, and we have
$$
\frac{d\hat\ce_{n,k,t}}{dt}=
\frac{1}{2t^{\frac{1}{2}}}
\sum_{l=1}^n\Big(  \beta u_N\ q^{\frac{p-1}{2}}
\sum_{J \in \bar Q_{N,k}^p} z_{J}\ \ep^l_{J}\Big)\
\hat\ce_{n,k,t}.
$$
Then,
\begin{eqnarray*}
\hat\nu_{k,t}'(f)=\frac{\beta u_N\ q^{\frac{p-1}{2}}}{2t^{\frac{1}{2}}}
&\be\Bigg[& \hat Z_{(k),t}^{-n}
\sum_{l=1}^n\av \Big(f 
\sum_{J \in \bar Q_{N,k}^p} z_{J}\ \ep^l_{J}\ \hat\ce_{n,k,t}\Big)
\\
&&-n \hat Z_{(k),t}^{-(n+1)}
\av \Big(f 
\sum_{J \in \bar Q_{N,k}^p} z_{J}\ \ep^{n+1}_{J}\ \hat\ce_{n+1,k,t}\Big)
\Bigg].
\end{eqnarray*}
An integration by parts formula with respect to the random variable $z_{J}$ implies
\begin{multline*}
\hat\nu_{k,t}'(f)=\frac{\beta^2 u_N^2 q^{p-1}}{2}
\sum_{J \in \bar Q_{N,k}^p}\Bigg[
\sum_{l=1}^n\sum_{l'=1}^n \nu_{k,t}\lp f \ep_{J}^{l} \ep_{J}^{l'}  \rp
-\sln \nu_{k,t}\lp f \ep_{J}^{l} \ep_{J}^{n+1}  \rp\\
-\sum_{l=1}^{n+1} \nu_{k,t}\lp f \ep_{J}^{l} \ep_{J}^{n+1}  \rp
+n(n+1)\nu_{k,t}\lp f \ep_{J}^{n+1} \ep_{J}^{n+2}  \rp \Bigg]
\end{multline*}
Now, from Lemma \ref{cardbn}, we obtain easily 
$$u_N^2 |\bar Q_{N,k}^p| \le \frac{K}{N},$$
that gives us  the desired result.
\hfill
$\Box$

%\vspace{0.2cm}
 
\noindent
{\bf Proof of Lemma \ref{l2.5}:}
Since $f^-$ depends on $\{\si_1^l,\ldots,\si_{N-p}^{l};l$ $\le n  \}$, invoking Lemma \ref{l2.4}, we only
need to work with 
$\hat\nu_{p-1,0}\big(\prod_{l\le n, m_l\in \mathcal{M}_l}
\ep_{m_l}^l\big)$.
Indeed,
$$
\nu_{p-1,0}\Big(f^-\prod_{l\le n,  m_l\in \mathcal{M}_l}
\ep_{m_l}^l\Big)= \nu_{p-1,0}(f^-)
\ \nu_{p-1,0}\Big(\prod_{l\le n, m_l\in \mathcal{M}_l}
\ep_{m_l}^l\Big),
$$
and 
$$
\left|\nu_{p-1,0}\Big(\prod_{l\le n, m_l\in \mathcal{M}_l}
\ep_{m_l}^l\Big)
-\hat\nu_{p-1,0}\Big(\prod_{l\le n, m_l\in \mathcal{M}_l}
\ep_{m_l}^l\Big)  \right|
\le \frac{K(\beta)}{N}.
$$
We will now divide our proof in two steps:

%\vspace{0.2cm}

\noindent
{\bf Step 1:} By the construction of $\hat\nu_{p-1,0}$, 
$\hat \ce_{n,p-1,0}$ and $\mathcal{M}_l$, using independence (of the $\ep_j$ with respect to the uniform measure on $\{-1;1\}^{np}$ and of the random variables $\hat z_j$), we get
\begin{align*}
&\hat\nu_{p-1,0}\Big(
\prod_{l\le n , m_l\in \mathcal{M}_l}
\ep_{m_l}^l
\Big)=\be\Bigg[\frac{\av\big(
\prod_{l\le n , m_l\in \mathcal{M}_l}
\ep_{m_l}^l\ \hat \ce_{n,p-1,0}\big)}{\hat Z_{(p-1)}^n}\Bigg]\\
&\qquad =\be\Bigg[\prod_{l\le n , m_l\in \mathcal{M}_l}  
\tanh \big(\beta u_N q^{\frac{p-1}{2}} \hat z_{m_l}+h\big)\Bigg]\\
&\qquad = \be\Bigg[\prod_{j\le p-1}\Big\{
\tanh \big(\beta u_N q^{\frac{p-1}{2}} \hat z_{j}+h\big)
\Big\}^{\sum_{l=1}^n \mathbb{I}_{\{j\in \mathcal{M}_l\}}}\Bigg]\\
&\qquad = \prod_{j\le p-1} \be \Bigg[
\tanh^{\sum_{l=1}^n \mathbb{I}_{\{j\in \mathcal{M}_l\}}}
\Big(\beta u_N
q^{\frac{p-1}{2}} \hat z_j+h \Big)  \Bigg].
\end{align*}

%\vspace{0.2cm}

\noindent
{\bf Step 2:}
By Lemma \ref{cardbn} and the fact that $\hat z_j$ is a centered Gaussian random
variable with variance ${N-k \choose p-1}$, we have
\begin{equation}\label{estep2}
E\lc u_N^2 \hat z_j^2\rc=\frac{p}{2}+\mathcal{O}\Big(\frac{1}{N}\Big).
\end{equation}
For $s>0$, set now $\psi(s)=\be[\tanh^m(X_s+h)]$, where $X_s$ is a centered Gaussian random variable
with variance $s^2$. Then
\begin{eqnarray*}
\psi(s)&=&\frac{1}{\sqrt{2\pi s^2}} \int_{-\infty}^\infty
\tanh^m(u+h)\ e^{-\frac{1}{2s^2}u^2} \ du\\
&=&\frac{1}{\sqrt{2\pi}} \int_{-\infty}^\infty
\tanh^m(vs+h)\ e^{-\frac{1}{2}v^2} \ dv.
\end{eqnarray*}
Note that $|\psi'(s)|\le K(m)$. Then, using  the fact that $Y$ 
is a standard Gaussian random variable and (\ref{estep2}), we have
\begin{eqnarray*}
&&\left| 
\be\lc \tanh^m\lp \beta\lp \frac{p}{2} \rp^{\frac{1}{2}}q^{\frac{p-1}{2}}\ Y+h \rp  \rc
-\be\lc \tanh^m\lp \beta u_N\ q^{\frac{p-1}{2}}\ \hat z_j+h  \rp \rc \ \right|\\
&&\quad\quad\quad = \left|\ \psi\Big(\beta q^{\frac{p-1}{2}}  
\sqrt{\frac{p}{2}}\Big) - 
\psi\Big(\beta q^{\frac{p-1}{2}}  u_N \sqrt{\be(\hat z_j^2)}\Big) \ \right|\\ 
&&\quad\quad\quad\le K(m) \beta  q^{\frac{p-1}{2}}
\left| \sqrt{\frac{p}{2}}- \sqrt{u_N^2\be(\hat z_j^2)} \right|\\
&&\quad\quad\quad\le  \frac{K(m,\beta)}{N},
\end{eqnarray*}
which shows our claim.
\hfill
$\Box$

\subsection{Particular case: The 1-cavity}\label{ucav}

The results of the previous subsection sometimes take a simpler form when expressed for a cavity of order 1. It is then useful to summarize them in this particular case:
let 
$$
\beta_{-}=\lp \frac{N-1}{N} \rp^{\frac{p-1}{2}}\beta,
$$
and $\langle \cdot \rangle_-$ the averaging with respect to the Gibbs
measure on $\Sigma_{N-1}$ at inverse temperature $\beta_{-}$.

Let $n\ge 1$ and $\si^1,\dots,\si^n$ be $n$ independent copies of a $N$-spins
configuration. For any $j\in \{1,\dots,n\}$, we denote $\sigma^j=(\rho^j,\ep^j)$, where
$\rho^j \in  \Sigma_{N-1}$ and $\ep^j\equiv\ep_1^j \in \{-1,1\}$. Set
$$Q_{N,1}^p= \big\{J=(i_1,\dots,i_{p-1},N)\in \mathbb{N}^p; 1 \le i_1<\cdots <i_{p-1}\le N-1\big\},$$
and, in this case, for $J\in Q_{N,1}^p$, 
$$
\eta_J=\si_{i_1}\cdots\si_{i_{p-1}},\qquad\quad
\ep_J=\si_{N}=\ep,
$$
and
\begin{displaymath}
g_{(1)}(\eta,\ep)
=\ep\ g(T(\rho)),
\end{displaymath}
being
$$g(T(\rho))=\beta u_N \sum_{J\in Q_{N,1}^p} g_J\ \eta_J.$$
Here, for one configuration, we have
\begin{displaymath}
-H_{N,\beta,h}(\si)=-H_{N-1,\beta_-,h}(\rho) + 
\ep\big[g(T(\rho)) + h\big].
\end{displaymath}
And we can also define
\begin{displaymath}
g_{(1),t}(\eta,\ep)= \ep \ g_t(T(\rho)),
\end{displaymath}
with
\begin{displaymath}
 g_t(T(\rho))=
 t^{\frac{1}{2}}g(T(\rho)) 
+\beta u_N q^{\frac{p-1}{2}}(1-t)^{\frac{1}{2}}
\sum_{J\in Q_{N,1}^p} z_J
\end{displaymath}
where $t\in[0,1]$, $q\in [0,1]$ and where $\{z_J ;J\in  Q_{N,1}^p \}$ 
is a family of independent standard Gaussian random variables, 
also independent of all the disorder $g$. 

Let us write
\begin{eqnarray*}
\ce_{n,1,t}&=&\exp\Big\{ \sum_{l=1}^{n}\ep^l 
\big[g_t(T(\rho^l))+h\big]\Big\},\\
Z_{(1),t}&=&\langle \av
\ce_{1,1,t} \rangle_{-}
=\langle \cosh \big[g_t(T(\rho^l))+h\big]\rangle_{-}. 
\end{eqnarray*}
For $f:\ssn^n\longrightarrow \mathbb{R}$, we can define
\begin{eqnarray*}
\langle f \rangle_{1,t}&=&\frac{\langle\av f \ce_{n,1,t} \rangle_{-}}{Z_{(1),t} ^n},\\
\nu_{1,t}(f)&=& \be\langle f \rangle_{1,t}.
\end{eqnarray*}
Then, for $t\in[0,1]$ and $f:\Sigma_N^n\to\R$, the derivative of $\nu_{1,t}(f)$ with respect to $t$ takes the exact form of relation (\ref{difnuk}) with $p=1$. Moreover,  as a particular case of Proposition \ref{posmaj}, we also get, for a non-negative function $f$ and $N$ large enough, that
\begin{equation}\label{efita}
\nu_{1,t}(f)\le \exp\big\{4\beta^2 n^2 p\big\}  \nu(f).
\end{equation}
An important remark is the fact that in order to prove the equivalent to
Lemma \ref{l2.5} (which will also be given in a simpler form), we do not need Lemma \ref{l2.4}:
let $f^-:\ssm^n\to\R$ be a function depending 
on $\{\si_1^l,\ldots,\si_{N-1}^{l};l$ $\le n  \}$. 
Let $Y$ be a standard Gaussian random variable. 
Then, for a constant $L_1(\beta)>0$,
\begin{equation}\label{e2.3.3}
\Big|\nu_{1,0}\big(f^- \ep^1 \cdots \ep^n\big)
-\nu_{1,0}(f^-) \be\Big[ \tanh^n
\big(\beta \lp\frac{p}{2} \rp^{\frac{1}{2}}
q^{\frac{p-1}{2}}Y+h \big)  \Big]\Big|\le \frac{L_1(\beta)}{N}\|f^-\|_\infty.
\end{equation}     
\section{Behavior of the overlap}

In this section we will study  the limiting behavior of the overlap
of two configurations, namely
$$R_{l,l'}=\frac{1}{N}\sum_{i=1}^N \si_i^l \si_i^{l'},$$
where $\si_i^l$ and $\si_i^{l'}$ are understood as two independent configurations under 
$G_N$. In the sequel, the following assumption on $\beta$, that determines our high temperature region, will have to be made:
\begin{itemize}
\item[{\bf (H)}]The parameter $\beta>0$ is smaller than a constant $\beta_p$ defined by
$$
8p^2\beta_p^2\exp\lp 16 \beta_p^2 p \rp=\frac12.
$$
\end{itemize}
We will see then that the constant $q=q_p$ which will be the $L^2$ limit of $R_{1,2}$, is the unique solution to the equation
\begin{equation}\label{qq}
q=\be \Bigg[ \tanh^2\Bigg(\beta \Big(\frac{p}{2}\Big)^{\frac{1}{2}}
q^{\frac{p-1}{2}} Y +h\Bigg)\Bigg].
\end{equation}
Observe that $\hat q_2=q$, where $\hat q_n$ is defined by relation (\ref{defqn}).

First, at Subsection \ref{selfav}, we will obtain the self averaging result for $R_{1,2}$, one of the main steps 
towards the replica symmetric formula. Then, at Subsection \ref{expmom}, using an elaboration of the arguments of Subsection \ref{selfav}, we will get the existence of exponential moments for $N(R_{1,2}-q)$. This will allow us, at Subsection \ref{lurate}, using higher order expansions, to get a sharp bound for the quantity $\nu(R_{1,2})-q$.

Along this paper we will use the following two deterministic results on the overlap. The first one 
is taken from Talagrand \cite[Lemma 5.11]{Tp}:
\begin{equation}\label{e3.1}
\left| u_N^2\sum_{J\in Q_{N,1}^p}\eti^l\eti^{l'}
-\frac{p}{2}R_{l,l'}^{p-1} \right|\le\frac{K}{N},
\end{equation}
while the second one is an easy consequence of  Lemma \ref{cardbn}:
\begin{equation}\label{e3.2}
\left| u_N^2|Q_{N,1}^p|q^{p-1}-\frac{p}{2}q^{p-1} \right|
\le
\frac{K}{N}.
\end{equation}

Let us first state an elementary Proposition, that will give some useful information about the whole $p$ spins system.
\begin{proposition}
Under assumption {\rm (H)}, equation (\ref{qq}) has a unique solution $q_p$ in $[0,1]$. Moreover
$$
\lim_{p\to\infty}q_p=\tanh^2(h).
$$
\end{proposition}

%\vspace{0.2cm}

\noindent
{\bf Proof:}
Let $\vp,\phi_p:\ou\to\ou$ be defined by
$$
\vp(x)=x,\quad \phi_p(x)=\be \Bigg[ \tanh^2\Bigg(\beta \Big(\frac{p}{2}\Big)^{\frac{1}{2}}
x^{\frac{p-1}{2}} Y +h\Bigg)\Bigg].
$$
It is easily seen that $\vp(0)=0$, $\phi_p(0)=\tanh^2(h)$ on one hand, and that $\vp(1)=1$, $\phi_p(1)<1$ on the other hand. Furthermore, a simple Gaussian integration by parts argument (see also the proof of Lemma \ref{l4.2}) shows that
$$
\phi^{'}_p(x)=\frac{\beta^2p(p-1)}{2}x^{p-2}
\be\lc \psi\lp \beta \lp\frac{p}{2}\rp^{1/2}x^{\frac{p-1}{2}}Y+h \rp \rc,
$$
where
$$
\psi(u)=\frac{1-2\sinh^2(u)}{\cosh^4(u)}.
$$
A quick study of $\psi$ shows that $\|\psi\|_\infty=1$, and hence, for any $x\in\ou$,
$$
|\phi'_p(x)|\le \frac{\beta^2p(p-1)}{2}x^{p-2}
\le 8p^2\beta_p^2\exp\lp 16 \beta_p^2 p \rp.
$$
Hence, if $\beta$ satisfies condition (H), the existence and uniqueness of the solution to (\ref{qq}) is trivially obtained. The second claim easily follows from the fact that $\phi_p(0)=\tanh^2(h)$, and that for any $a\in (0,1)$
$$
\lim_{p\to\infty}\sup_{x\in [0,a)}|\phi^{'}_p(x)|=0.
$$
\hfill
$\Box$

\subsection{Self averaging property}\label{selfav}
This part of the paper is devoted to prove that $R_{1,2}$ converges to $q$ in a $L^2$ sense. A 1-cavity will be enough to reach the conclusion of this Section, and we refer to Section \ref{ucav} for further notations and results on this method.
 First of all recall that
$\nu_{1,1}(f)=\nu(f)$. 

\begin{proposition}\label{p3.1}
Let $f$ be a function from 
$\ssn^n$ to $\R$, and $\al_1,\al_2>1$ such that 
$\al_1^{-1}+\al_2^{-1}=1$. Then, 
there exists a positive constant $L_2$ such that
\begin{multline*}
\left| \nu(f)-\nu_{1,0}(f) \right|\le
(n p\beta)^2\exp\lp 4\beta^2n^2p \rp\\
\nu^{1/\al_1}(|f|^{\al_1})
\lp\nu^{1/\al_2}(|R_{1,2}-q|^{\al_2})+\frac{L_2}{N}  \rp.
\end{multline*}
\end{proposition}

%\vspace{0.2cm}

\noindent
{\bf Proof:}
Consider the term
$$
U_t=\beta^2u_N^2\sum_{J\in Q_{N,1}^p}\
\nu_{1,t}\lp f\ep^{l}\ep^{l'}(\eti^{l}\eti^{l'}-q^{p-1})\rp.
$$
Then, by H\"older's inequality,
$$
U_t\le\beta^2\nu_{1,t}^{1/\al_1}(|f|^{\al_1})
\nu_{1,t}^{1/\al_2}\lp \Big| u_N^2\sum_{J\in Q_{N,1}^p}\ \eti^l\eti^{l'}
-u_N^2|Q_{N,1}^p|q^{p-1} \Big|^{\al_2} \rp.
$$
Using (\ref{e3.1}) and (\ref{e3.2}), we obtain
$$
U_t\le\frac{p\beta^2}{2}\nu_{1,t}^{1/\al_1}(|f|^{\al_1})
\lp \nu_{1,t}^{1/\al_2}\lp \left| R^{p-1}_{l,l'}-q^{p-1} \right|^{\al_2}\rp + \frac{L_2}{N} \rp,
$$
and since $|R^{p-1}_{l,l'}-q^{p-1}|\le p|R_{l,l'}-q|$, we get
$$
U_t\le\frac{p^2\beta^2}{2}\nu_{1,t}^{1/\al_1}(|f|^{\al_1})
\lp \nu_{1,t}^{1/\al_2}\lp \left| R_{1,2}-q \right|^{\al_2}\rp + \frac{L_2}{N} \rp.
$$
By relation (\ref{majnutnu}), we then have
$$
U_t\le\frac{p^2\beta^2}{2}\exp\lp 4\beta^2n^2p \rp
\nu^{1/\al_1}(|f|^{\al_1}) \lp \nu^{1/\al_2}\lp \left| R_{1,2}-q \right|^{\al_2}\rp + \frac{L_2}{N}\rp.
$$
Our result is then obtained by iteration of this kind of calculations for the other terms in (\ref{derivnu}).

\hfill
$\Box$
%\vspace{0.1cm}

\begin{proposition}\label{p3.2} Let q be the solution to (\ref{qq}).
If $\beta$ satisfies {\rm (H)}, then
%There exists $\beta_p>0$ such that, for any $\beta\le \beta_p$, we have
$$\nu\big((R_{1,2}-q)^2\big)=\be \langle (R_{1,2}-q)^2\rangle \le \frac{K}{N}.$$
\end{proposition}
\vspace{0.2cm}

\noindent
{\bf Proof:} The symmetry between sites implies that
\begin{equation}\label{e3.6}
\nu\big((R_{1,2}-q)^2\big)=
\nu(\bar f),
\end{equation}
where 
$$\bar f= (\ep^1 \ep^2-q)(R_{1,2}-q)=A_1+A_2,$$
with
\begin{eqnarray*}
A_1&=&\frac{1}{N}(\ep^1 \ep^2-q)^2,\\ 
A_2&=& (\ep^1 \ep^2-q)\Big(R_{1,2}^--\frac{N-1}{N}q\Big),
\end{eqnarray*}
and 
$$R_{1,2}^-=\frac{1}{N}\sum_{i=1}^{N-1} \si_i^1 \si_i^{2}.$$
Since $|\ep^1 \ep^2-q|\le 2$, it is obvious that
\begin{displaymath}
\nu_{1,0}(A_1)\le \frac{4}{N}.
\end{displaymath}
On the other hand, by relation (\ref{e2.3.3}) and the fact that $q$ is the solution to (\ref{qq}), we get
\begin{eqnarray*}
\nu_{1,0}(A_2)&=&\nu_{1,0}\Big(R_{1,2}^--\frac{N-1}{N}\ q\Big)\\
&&\quad \times
\Bigg[\be\Big[ \tanh^2
\big(\beta \lp\frac{p}{2} \rp^{\frac{1}{2}}
q^{\frac{p-1}{2}}Y+h \big)  \Big]-\ q\Bigg]+\mathcal{O}\Big(\frac{1}{N}\Big)\\
&=&\mathcal{O}\Big(\frac{1}{N}\Big),
\end{eqnarray*}
and Proposition \ref{p3.1} for $n=\alpha_1=\alpha_2=2$ yields
\begin{displaymath}
\big|\nu(\bar f)-\nu_{1,0}(\bar f)\big|\le(2p\beta)^2
\exp\{16\beta^2p\}\ \nu^{1/2}(|\bar f|^{2})
\lp\nu^{1/2}(|R_{1,2}-q|^{2})+\frac{L_2}{N}  \rp.
\end{displaymath}
Then (\ref{e3.6}) and the estimates
for $A_1$ and $A_2$ imply
$$\nu\big((R_{1,2}-q)^2\big)\le 8p^2\beta^2 \exp\{16\beta^2p\}\
\nu\big((R_{1,2}-q)^2\big)+\frac{K}{N}.$$
Thus, if $\beta_p$ satisfies (H),
%is choosen so that
%$$8p^2\beta_p^2 \exp\{16\beta_p^2p\}\le \frac{1}{2},$$
we obtain the desired inequality.

\hfill
$\Box$
%\vspace{0.1cm}

\subsection{Exponential moments}\label{expmom}
The aim of this subsection is to bound the higher moments of $R_{1,2}-q$.
Notice that these bounds will be used in the next subsection in order to control 
$\nu(R_{1,2}-q)$.
\begin{theorem}\label{t3.3}
Let q be the solution to (\ref{qq}).
%There exists $\beta_p>0$ such that, for any $\beta\le \beta_p$
If $\beta$ satisfies {\rm (H)}, we have
$$
%\begin{equation}\label{e3.2.1}
\nu\big((R_{1,2}-q)^{2l}\big)
=\be \langle (R_{1,2}-q)^{2l}\rangle \le \Big(\frac{Ll}{N}\Big)^l,
$$
%\end{equation}
where $L$ does not depend on $l$.
\end{theorem}
%\vspace{0.1cm}
\begin{remark}
Theorem \ref{t3.3} implies that there exists $M>0$ such that
$$\nu\Bigg(\exp\Big\{\frac{N}{M}(R_{1,2}-q)^2\Big\}
\Bigg)\le M,$$
and hence the title of this section.
Indeed, this is an immediate consequence of the equality 
$e^{x^2}=\sum_{l\ge 0} \frac{x^{2l}}{l!}$
and the fact $\Big(\frac{l}{3}\Big)^l\le l!\le l^l$. 
\end{remark}

The proof of Theorem \ref{t3.3} goes along the same lines as
Theorem 2.5.1 in \cite{Tbk}, except for the introduction of a two steps induction due to the high number of interactions between spins. We will try to stress mainly on this difference.
We will proceed by induction over $l$, and the induction hypothesis 
will be
\begin{equation}\label{e3.2.2}
\nu\big((R_{1,2}-q)^{2\bar l}\big)
 \le \Big(\frac{L_0\bar l}{N}\Big)^{\bar l},\qquad \textrm{for any}\ \bar l
 \in\{1,\dots,l\},\end{equation}  
 being $L_0$ a fixed number.  The case $l=1$ has been proved in 
 Proposition \ref{p3.2}, and
if $L_0$ is large enough, we will show that
\begin{equation}\label{e3.2.3}
\nu\big((R_{1,2}-q)^{2l+2}\big)
 \le \Big(\frac{L_0(l+1)}{N}\Big)^{l+1}.\end{equation}
First of all, since $|R_{1,2}-q|\le 2$, for any
$\bar l\ge N$ assuming $L_0\ge 4$, we have
$$\nu\big((R_{1,2}-q)^{2\bar l}\big)\le 4^{\bar l}
 \le L_0^{\bar l}\Big(\frac{\bar l}{N}\Big)^{\bar l}.$$
So, we can suppose $l\le N-1$.

In order to prove (\ref{e3.2.3}) we will need the following lemma.

\begin{lemma}\label{lemind}
Assume (\ref{e3.2.2}) and $\ l\le N-1$. Then, if $L_0\ge 4$ we
have
\begin{eqnarray*}
\nu\big(|R_{1,2}-q|^j\big)&\le& \Big(\frac{L_0(j+1)}{N}\Big)^{j/2},
\qquad \forall j\le 2l,\\
\nu\big((R_{1,2}^--q)^{2l}\big)&\le& 
3\Big(\frac{L_0(l+1)}{N}\Big)^{l}.
\end{eqnarray*}
\end{lemma}
%\vspace{0.1cm}

\noindent {\bf Proof:} See the proof of Lemma 2.5.1 in \cite{Tbk}.
\hfill
$\Box$

%\vspace{0.1cm}

\noindent {\bf Proof of Theorem \ref{t3.3}:}
Our goal is to prove (\ref{e3.2.3})
assuming (\ref{e3.2.2}). 
By symmetry we have
\begin{equation}\label{e3.2.6}
\nu\big((R_{1,2}-q)^{2l+2}\big)  
=\nu(\bar f)=\nu_{1,0}(\bar f)+\big[\nu(\bar f)-\nu_{1,0}(\bar f)\big],
\end{equation}
where
$$\bar f=(\ep^1 \ep^2 -q)\big(R_{1,2}-q\big)^{2l+1}.$$
Applying Proposition \ref{p3.1} 
with $n=2$, $\alpha_1=\frac{2l+2}{2l+1}$ and $\alpha_2=2l+2$,
and using $|\ep^1 \ep^2-q|\le 2$  we obtain
\begin{equation}\label{eem1}
\big|\nu(\bar f)-\nu_{1,0}(\bar f)\big|
\le 8p^2\beta^2 e^{16\beta^2p}\Big[\nu\big((R_{1,2}-q)^{2l+2}\big)
+ \frac{L_2}{N}\nu^{\frac{2l+1}{2l+2}}\big((R_{1,2}-q)^{2l+2}\big)\Big].
\end{equation}
Assuming condition (H) and plugging (\ref{eem1}) into (\ref{e3.2.6}) we get,
for $\beta \le \beta_p$,
\begin{equation}\label{ee19}
\nu\big((R_{1,2}-q)^{2l+2}\big)
\le 2 \nu_{1,0}(\bar f)
+  \frac{L_2}{N}
\nu^{\frac{2l+1}{2l+2}}\big((R_{1,2}-q)^{2l+2}\big).
\end{equation}
This inequality, which was sufficient in the case of the SK model (see \cite{Tbk}), does not allow us to reach our conclusion here, and we will have to perform a second step in our induction: using that
\begin{displaymath}
(x+y)^\alpha\le x^\alpha +  y^\alpha,
\end{displaymath}
for $x,y\ge 0$ and $\alpha\in(0,1)$, 
from (\ref{ee19}) we easily obtain, for $\beta \le \beta_p$,
\begin{equation}\label{ee19bis}
\nu\big((R_{1,2}-q)^{2l+2}\big) 
\le A_1+A_2+A_3,
\end{equation}
where
\begin{eqnarray*}
A_1&=&2 \nu_{1,0}(\bar f),\\
A_2&=& \frac{L_2}{N} \ 2^{\frac{2l+1}{2l+2}} \
\nu_{1,0}^{\frac{2l+1}{2l+2}}(\bar f),\\
A_3&=& \frac{L_2}{N} \Bigg( \frac{L_2}{N} \
\nu^{\frac{2l+1}{2l+2}}\big((R_{1,2}-q)^{2l+2}\big)
\Bigg)^{\frac{2l+1}{2l+2}}.
\end{eqnarray*}
Let us study first $A_3$.
Using (\ref{e3.2.2}) we get for  $\beta \le \beta_p$ 
\begin{eqnarray}
A_3&\le&\frac{4\ L_2^2 \
\Big[\nu\big((R_{1,2}-q)^{2l}\big)\Big]^{\left(\frac{2l+1}{2l+2}\right)^2}}
{N^{1+\frac{2l+1}{2l+2}}}\nonumber\\
&\le& \frac{4\ L_2^2\ (L_0\ l)^l}{N^{1+\frac{2l+1}{2l+2}+
l\left(\frac{2l+1}{2l+2}\right)^2}}\nonumber\\
&\le& \frac{4\ L_2^2\ (L_0\ l)^l}{N^{l+1}}.\label{ee20}
\end{eqnarray}
In order to study $A_1$ note that
\begin{eqnarray*}
|\nu_{1,0}(\bar f)|&\le&
\Big|\nu_{1,0}\Big((\ep^1 \ep^2 -q)
\big[(R_{1,2}-q)^{2l+1}-(R_{1,2}^--q)^{2l+1}\big]\Big)\Big|\\
&&+ \Big|\nu_{1,0}\big((\ep^1 \ep^2 -q)(R_{1,2}^--q)^{2l+1}\big)\Big|.
\end{eqnarray*}
The independence between $\ep^1 \ep^2$
and $R_{1,2}^-$ under $\nu_{1,0}$, inequalities (\ref{e2.3.3}) and (\ref{efita}) and  Lemma \ref{lemind} yield 
\begin{eqnarray*}
\Big|\nu_{1,0}\big((\ep^1 \ep^2 -q)(R_{1,2}^--q)^{2l+1}\big)\Big|
&\le& \Big|\nu_{1,0}(\ep^1 \ep^2 -q)\Big|
\Big|\nu_{1,0}\big((R_{1,2}^--q)^{2l+1}\big)\Big|\nonumber\\
&\le& \frac{L_1\ e^{16 \beta_p^2 p}}{2N}\
\big|\nu\big((R_{1,2}^--q)^{2l+1}\big)\big| 
\nonumber\\
&\le& 3\ L_1\ e^{16 \beta_p^2 p}\ \frac{(L_0(l+1))^l}{N^{l+1}},
\end{eqnarray*}
with $L_1:=L_1(\beta_p)$ given by (\ref{e2.3.3}). 
On the other hand, using the inequality
$|x^{2l+1}-y^{2l+1}|\le (2l+1)|x-y|(x^{2l}+y^{2l})$ 
%Lemma \ref{lemind},
and similar arguments as before,
we get  

\begin{eqnarray*}
&&\Big|\nu_{1,0}\Big((\ep^1 \ep^2 -q)
\big[(R_{1,2}-q)^{2l+1}-(R_{1,2}^--q)^{2l+1}\big]\Big)\Big|\nonumber\\
&&\qquad\qquad\qquad\le 
2
\Big|\nu_{1,0}\Big((R_{1,2}-q)^{2l+1}-(R_{1,2}^--q)^{2l+1}\Big)\Big|\nonumber\\
&&\qquad\qquad\qquad\le 
\frac{2(2l+1)}{N}
\Big[\nu_{1,0}\Big((R_{1,2}-q)^{2l}\Big)
+\nu_{1,0}\Big((R_{1,2}^--q)^{2l}\Big)\Big]\nonumber\\
&&\qquad\qquad\qquad\le
\frac{2(2l+1)\ e^{16\beta_p^2p}}{N} 
\Big[\nu\Big((R_{1,2}-q)^{2l}\Big)+
\nu\Big((R_{1,2}^--q)^{2l}\Big)\Big]\nonumber\\
&&\qquad\qquad\qquad\le
\frac{2(2l+1)\ e^{16\beta_p^2p}}{N} 
\Big[\Big(\frac{L_0 l}{N}\Big)^l +
3\Big(\frac{L_0 (l+1)}{N}\Big)^l\Big]\\
&&\qquad\qquad\qquad
\le\frac{4\ e^{16\beta_p^2p}\ 
\big[(L_0 l)^l + 3 (L_0(l+1))^l\big](l+1)}{N^{l+1}}\\
&&\qquad\qquad\qquad\le
8\ e^{16\beta_p^2 p}\ \frac{L_0^l(l+1)^{l+1}}{N^{l+1}}. 
\end{eqnarray*}   
So, it follows that
\begin{eqnarray}\label{ee21}
A_1&\le& 2\ e^{16\beta_p^2p}\ (3 L_1+8)\
\frac{L_0^l (l+1)^{l+1}}{N^{l+1}}\nonumber\\
&\le& 
2\ e^{16p}\ (3 L_1+8)\
\frac{L_0^l (l+1)^{l+1}}{N^{l+1}}.
\end{eqnarray}
It is also easy to check that
\begin{eqnarray}\label{ee22}
A_2&\le &
2\ \frac{L_2}{N} \Bigg[2\ e^{16p}\ (3 L_1+8)
\frac{L_0^l\ (l+1)^{l+1}}{N^{l+1}}\Bigg]^{\frac{2l+1}{2l+2}}\nonumber\\
&\le&\frac{4\ L_2\ e^{16p}\ (3 L_1+8)\ L_0^l\ (l+1)^{l+1}}
{N^{1+ (l+1)\left(\frac{2l+1}{2l+2}\right)}}\nonumber\\
&\le&\frac{4\ L_2 \ e^{16p}\ (3 L_1+8)\ L_0^l\ (l+1)^{l+1}}
{N^{l+1}}.
\end{eqnarray} 
Puting together (\ref{ee19bis}), (\ref{ee20}), (\ref{ee21}) and (\ref{ee22})
we obtain, for $\beta \le \beta_p$,
\begin{displaymath}
\nu\big((R_{1,2}-q)^{2l+2}\big)\le \bar K
\frac{ L_0^l\ (l+1)^{l+1}}{N^{l+1}},
\end{displaymath}
with 
$$\bar K=4\ L_2^2\ +\ 2\ e^{16p}\ (3L_1+8)
\ (1+2 L_2).$$ 
So, if $L_0\ge \bar K$ the proof is completed. 
\hfill
$\Box$

%\vspace{0.1cm}

\noindent
An easy consequence of Theorem \ref{t3.3} is the following

\begin{corollary}\label{c3.5}
Let q be the solution of  (\ref{qq}).
If  $\beta$ satisfies {\rm (H)}, we have
$$
%\begin{equation}\label{e3.2.7}
\be \langle (R_{1,2}^{p-1}-q^{p-1})^{2l}\rangle 
\le \Big(\frac{\hat Ll}{N}\Big)^l.
$$
%\end{equation}
\end{corollary}
%\vspace{0.1cm}
\noindent {\bf Proof:} This result is an obvious consequence of 
Theorem \ref{t3.3}. Indeed,
\begin{eqnarray*}
\be \langle (R_{1,2}^{p-1}-q^{p-1})^{2l}\rangle  
&\le& \be \Big\langle (R_{1,2}-q)^{2l}
\big[\sum_{k=0}^{p-2}R_{1,2}^{p-2-k}q^k\big]^{2l}\Big\rangle\\
&\le& (p-1)^{2l}\ \be \langle (R_{1,2}-q)^{2l}\rangle\\
&\le& \Big(\frac{L(p-1)^2 l}{N}\Big)^l,
\end{eqnarray*}
which is the desired result.
\hfill
$\Box$

%\vspace{0.1cm}

\noindent
Another immediate Corollary of the expoenential inequalities for  the overlap is a useful result on the expansions of $\nu(f)$, that we label for further use.
\begin{corollary}\label{c3.6}
If  $\beta$ satisfies {\rm (H)}  we have, for a function $f$ on $\Sigma^n_N$ and
$k\ge 1$,
\begin{eqnarray}
&&\big|\nu(f)- \nu_{k,0}(f) \big|\le\frac{K}{N^{\frac{1}{2}}}
\nu^{\frac{1}{2}}(f^2),\label{e3.5}\\
&&\big|\nu(f)- \nu_{k,0}(f) - \nu'_{k,0}(f) \big|
\le\frac{K}{N}
\nu^{\frac{1}{2}}(f^2).\label{e3.6'}
\end{eqnarray}
\end{corollary}

%\vspace{0.1cm}            

\noindent {\bf Proof:} We refer to \cite{Tbk} for the proof of this corollary.
\hfill
$\Box$

\subsection{Upper bound for $\nu(R_{1,2}-q)$}\label{lurate}
The main goal of this part of the paper will be to prove the following Theorem, that gives a sharp rate of convergence of $R_{1,2}$ towards $q$.

\begin{theorem}\label{t3.8}
Let q be the solution of  (\ref{qq}).
%There exists $\beta_p>0$ such that, for any $\beta\le \beta_p$, we have
Then, if $\beta$ satisfies {\rm (H)}, we have
$$
\big|\nu\big(R_{1,2}-q\big)\big|
\le \frac{K}{N}.
$$
\end{theorem}
%\vspace{0.1cm}

\noindent
An immediate consequence of this theorem is the following result.

\begin{corollary}\label{cmajoration} 
Let q be the solution of  (\ref{qq}).
If $\beta$ satisfies {\rm (H)}, we have
$$
\big|\nu\big(R_{1,2}^m-q^m\big)\big|
\le \frac{K(m)}{N},
$$
for all $m\ge 1$.
\end{corollary}
%\vspace{0.1cm}
\noindent {\bf Proof:} 
For a fixed $m\ge 1$, by Taylor's expansion, we have
\begin{equation}\label{emajoration}
R_{1,2}^m=q^m+m\ q^{m-1}(R_{1,2}-q)+\frac{m(m-1)}{2}\ \xi^{m-2} (R_{1,2}-q)^2,
\end{equation}
where $\xi \in  (R_{1,2} \wedge q, R_{1,2} \vee q)$.  So
\begin{displaymath}
\nu(R_{1,2}^m-q^m)=m\ q^{m-1}\nu(R_{1,2}-q)+\frac{m(m-1)}{2}\ 
\nu\big(\xi^{m-2}  (R_{1,2}-q)^2\big).
\end{displaymath}
Since $|\xi|\le 1$ and  using Theorem \ref{t3.8} and Proposition \ref{p3.2}
we obtain
$$
\big|\nu\big(R_{1,2}^m-q^m\big)\big|
\le \frac{K(m)}{N}.
$$
\hfill
$\Box$

%\vspace{0.1cm}

\noindent
Conversely, the next proposition will show that, in order to prove Theorem \ref{t3.8}, it will be enough to establish the following upper bound:
\begin{equation}\label{e3.3.1}
\nu(R_{1,2}^{p-1}-q^{p-1})\le \frac{K}{N}.
\end{equation}
\begin{proposition}\label{p3.3.1}
For $N$ large enough, there exist positive constants $L_3$ and $L_4$
such that, for any $m\ge 1$,
$$|\nu(R_{1,2}-q)|\le L_3 |\nu(R_{1,2}^{m}-q^{m})|
+\frac{L_4}{N}.$$
\end{proposition}

%\vspace{0.1cm}

\noindent {\bf Proof:} 
By (\ref{emajoration}),
and since $q$ is a strictly positive number, we have
\begin{displaymath}
(R_{1,2}-q)=\frac{1}{m q^{m-1}}(R_{1,2}^m-q^m)
-\frac{m-1}{2q^{m-1}}\ \xi^{m-2} (R_{1,2}-q)^2,
\end{displaymath}
where  $\xi \in  (R_{1,2} \wedge q, R_{1,2} \vee q)$.
Using Proposition \ref{p3.2} we can bound $\nu[(R_{1,2}-q)^2]$ by $\frac{L}{N}$, finishing the proof.
\hfill
$\Box$

%\vspace{0.1cm}

We will now prepare the proof of (\ref{e3.3.1}) by a series of lemmas, beginning with some deterministic estimates for the overlap.

\begin{lemma}\label{nouu} 
\begin{eqnarray}
R_{1,2}^{p-1}&=&\frac{(p-1)!}{N^{p-1}}
\sum_{\hat J \in A_N^{p-1}}
\eta_{\hat J}^1  \eta_{\hat J}^2 
+ \mathcal{O}\Big(\frac{1}{N}\Big)\label{eS12}\\ 
&=& \frac{(p-1)!}{N^{p-1}} \sum_{\hat J \in \asp} 
\eta_{\hat J}^1  \eta_{\hat J}^2 
+ \mathcal{O}\Big(\frac{1}{N}\Big),\label{equasif}
\end{eqnarray}
\end{lemma}

%\vspace{0.1cm}
\noindent {\bf Proof:}
Let $N_{p-1}=
\{(i_1,\dots,i_{p-1})\in \{1,\dots,N\}^{p-1}\}$.
We can easily check
\begin{eqnarray}
R_{1,2}^{p-1}&=&\frac{1}{N^{p-1}}
\sum_{(i_1,\dots,i_{p-1})\in N_{p-1}}
\sigma^1_{i_1}\cdots\sigma^1_{i_{p-1}}
\sigma^2_{i_1}\cdots\sigma^2_{i_{p-1}}\label{etll}\\ 
&=&\frac{1}{N^{p-1}}\Bigg[\sum_{(i_1,\dots,i_{p-1})\in \bar N_{p-1}}
\sigma^1_{i_1}\cdots\sigma^1_{i_{p-1}}
\sigma^2_{i_1}\cdots\sigma^2_{i_{p-1}}\nonumber\\
&& 
+\sum_{(i_1,\dots,i_{p-1})\in \bar N^c_{p-1}} 
\sigma^1_{i_1}\cdots\sigma^1_{i_{p-1}}
\sigma^2_{i_1}\cdots\sigma^2_{i_{p-1}}\Bigg],\nonumber
\end{eqnarray}
where  $\bar N_{p-1}$ is the set of elements $(i_1,\dots,i_{p-1})$
belonging to $N_{p-1}$ such that all the elements $i_1,\dots,i_{p-1}$ are different
and $\bar N^c_{p-1}$ is the complementary set of $\bar N_{p-1}$
with respect to $N_{p-1}$, that means,
the set of elements $(i_1,\dots,i_{p-1})$ belonging to $N_{p-1}$ such that
at least two values $i_k, i_{k'},  
k \neq k'$, are equal. Then, Lemma \ref{cardnp} imply
(\ref{eS12}). Moreover, the equality 
$$
|A_N^{p-1} \cap (\asp)^c|=P_{p-2}(N)
$$
gives us (\ref{equasif}).
\hfill
$\Box$

\begin{corollary}\label{noudos} 
$$
%\begin{equation}\label{erelacio}
u_N^2 \sum_{\hat J\in \asp} 
\eta_{\hat J}^{1}\eta_{\hat J}^{2}
=\frac{p}{2} R_{1,2}^{p-1}+ \mathcal{O}\Big(\frac{1}{N}\Big).
$$
%\end{equation}
\end{corollary}

%\vspace{0.1cm}
\noindent {\bf Proof:}
Trivial from the definition of $u_N$.
\hfill
$\Box$

%\vspace{0.1cm}

\begin{lemma}\label{darrer} Let $f$ be a function from $\Sigma_N^n $ to $\R$, $k$ a positive integer, and $t\in\ou$.  Then 
\begin{eqnarray*}
& &\Big|u_N^2 \sum_{J\in Q_{N,p-1}^p}
\nu_{k,t} \Big( f(\eta_{J}^{1}\eta_J^{2}-q^{p-1})\
\ep_J^{1}\ep_J^{2}\Big)\\
& &\qquad - u_N^2 \sum_{l=1}^p   \sum_{{\J} \in \asp}
\nu_{k,t} \Big( f(\eta_{{\J}}^{1}\eta_{\J}^{2}-q^{p-1})\
\ep_l^{1}\ep_l^{2}\Big) \Big|
\le \frac{K}{N}(\nu_{k,t}(f^2))^{\frac{1}{2}}.\nonumber
\end{eqnarray*}
\end{lemma}

%\vspace{0.1cm}

\noindent {\bf Proof:}
Using Lemma \ref{cardbn} we have
$$u_N^2 |\bar Q_{N,p-1}^p| \le \frac{K}{N}.$$
Then, the decompostion of the set $Q_{N,p-1}^p$ given in (\ref{AA}) (see also Definition
\ref{setqmjr})  yields 
\begin{eqnarray*}
& &u_N^2 \sum_{J\in Q_{N,p-1}^p}
\nu_{k,t} \Big( f(\eta_{J}^{1}\eta_J^{2}-q^{p-1})\
\ep_J^{1}\ep_J^{2}\Big) \\
& & \qquad = u_N^2 \sum_{J\in \tilde Q_{N,p-1}^p}
\nu_{k,t} \Big( f(\eta_{J}^{1}\eta_J^{2}-q^{p-1})\
\ep_J^{1}\ep_J^{2}\Big) + \frac{K}{N}(\nu_{k,t}(f^2))^{\frac{1}{2}}.
\end{eqnarray*}
Note that we can write
$$
\tilde Q_{N,p-1}^p =\asp \times \{\ep_1,\dots,\ep_k\}.
$$
For $J \in  \tilde Q_{N,p-1}$, $\eta_J$ only depends on $\asp$, and $\ep_J$ is of the form $\ep_l$ for $l\le k$. So we can write  
$\eta_J=\eta_{\hat J}$ for $\hat J \in \asp$ 
instead of $\eta_J$  for $J \in  \tilde Q_{N,p-1}$.
Then, using
\begin{eqnarray*}
& &\sum_{J\in \tilde Q_{N,p-1}^p}
\nu_{k,t} \Big( f(\eta_{J}^{1}\eta_J^{2}-q^{p-1})\
\ep_J^{1}\ep_J^{2}\Big)\\
& &\qquad =
\sum_{l=1}^p   \sum_{{\J} \in \asp}
\nu_{k,t} \Big( f(\eta_{{\J}}^{1}\eta_{\J}^{2}-q^{p-1})\
\ep_l^{1}\ep_l^{2}\Big),
\end{eqnarray*}
the proof is completed.
\hfill
$\Box$

\begin{theorem}\label{t3.7}
Let q be the solution to (\ref{qq}).
If $\beta$ satisfies {\rm (H)},
%There exists $\beta_p>0$ such that, for any $\beta\le \beta_p$, 
we have
\begin{equation}\label{e3.3.7}
\big|\nu\big(R_{1,2}^{p-1}-q^{p-1}\big)\big|
\le \frac{K}{N}.
\end{equation}
\end{theorem}

%\vspace{0.1cm}

\noindent {\bf Proof:} Using (\ref{e3.1}), the symmetry among sites Lemma \ref{cardbn}  we get
\begin{eqnarray*}
\nu\big(R_{1,2}^{p-1}-q^{p-1}\big)&=&
\nu\Big(\frac{2}{p}u_N^2\sum_{J\in Q_{N,1}^p}\eta_J^1 \eta_J^2-q^{p-1}\Big)
+\mathcal{O}\Big(\frac{1}{N}\Big)\\
&=&\nu\Big(\frac{2}{p}u_N^2 |Q_{N,1}^p|  \ep_1^1\cdots
\ep_{p-1}^1 \ep_1^2\cdots \ep_{p-1}^2 - q^{p-1}\Big)
+\mathcal{O}\Big(\frac{1}{N}\Big)\\
&=&\nu(\bar f)+\mathcal{O}\Big(\frac{1}{N}\Big),
\end{eqnarray*}
where
$$\bar f=\prod_{j\le p-1} \ep^1_j \ep_j^2-q^{p-1}.$$
The estimate (\ref{e3.6'}) for $k=p-1$  yields
\begin{equation}
\big|\nu(\bar f)- \nu_{p-1,0}(\bar f) - \nu'_{p-1,0}(\bar f) \big|
\le\frac{K}{N}
\nu^{\frac{1}{2}}(\bar f^2).\label{e3.3.7'}
\end{equation}
Moreover, Lemma \ref{l2.5} and the  equation satisfied by $q$ imply
\begin{eqnarray}
\nu_{p-1,0}(\bar f)&=&\prod_{j\le p-1}
\be \Bigg[ \tanh^2
\Big(\beta \lp\frac{p}{2} \rp^{\frac{1}{2}}
q^{\frac{p-1}{2}}Y+h \Big)  \Bigg]-q^{p-1}
+\mathcal{O}\Big(\frac{1}{N}\Big)\nonumber\\
&=&q^{p-1}-q^{p-1}+  \mathcal{O}\Big(\frac{1}{N}\Big)\nonumber\\
&=&\mathcal{O}\Big(\frac{1}{N}\Big).\label{e3.3.10}
\end{eqnarray}

Let us now study $\nu'_{p-1,0}(\bar f)$. 
Applying (\ref{derivnu}) for
$k=p-1$, Lemma \ref{darrer} and the symmetry property among the $\ep_j$, we get 
\begin{eqnarray*}
\nu_{p-1,0}'(\bar f)&=&
\beta^2u_N^2\sum_{J\in Q_{N,p-1}^p}
\Bigg[\nu_{p-1,0} \Big( \bar f(\eta_{J}^{1}\eta_J^{2}-q^{p-1})\
\ep_J^{1}\ep_J^{2}\Big)\\
&&-4\nu_{p-1,0} \Big( \bar f(\eta_J^{1}\eta_J^{3}-q^{p-1})\
\ep_J^{1}\ep_J^{3}\Big)\\
&&+3
\nu_{p-1,0}\Big( \bar f(\eta_J^{3}\eta_J^{4}-q^{p-1})\ 
\ep_J^{3}\ep_J^{4}\Big)\Bigg]\\
&=&
\beta^2u_N^2\sum_{\hat J\in \asp}
\sum_{l=1}^{p-1}
\Bigg[\nu_{p-1,0} \Big( \bar f(\eta_{\hat J}^{1}\eta_{\hat J}^{2}-q^{p-1})\
\ep_{l}^{1}\ep_{l}^{2}\Big)\\
&&-4\nu_{p-1,0} \Big( \bar f(\eta_{\hat J}^{1}\eta_{\hat J}^{3}-q^{p-1})\
\ep_{l}^{1}\ep_{l}^{3}\Big)\\
&&+3
\nu_{p-1,0}\Big( \bar f(\eta_{\hat J}^{3}\eta_{\hat J}^{4}-q^{p-1})\ 
\ep_{l}^{3}\ep_{l}^{4}\Big)\Bigg]
+\mathcal{O}\Big(\frac{1}{N}\Big)\\
&=&
\beta^2u_N^2\sum_{\hat J\in \asp}
(p-1)\big[W_1-4W_2+3W_3\big]+\mathcal{O}\Big(\frac{1}{N}\Big),
\end{eqnarray*}
with
\begin{eqnarray*}
W_1&=& \nu_{p-1,0} \Big( \bar f(\eta_{\hat J}^{1}\eta_{\hat J}^{2}-q^{p-1})\
\ep_{p-1}^{1}\ep_{p-1}^{2}\Big),\\
W_2&=&
\nu_{p-1,0} \Big( \bar f(\eta_{\hat J}^{1}\eta_{\hat J}^{3}-q^{p-1})\
\ep_{p-1}^{1}\ep_{p-1}^{3}\Big),\\
W_3&=&
\nu_{p-1,0}\Big( \bar f(\eta_{\hat J}^{3}\eta_{\hat J}^{4}-q^{p-1})\ 
\ep_{p-1}^{3}\ep_{p-1}^{4}\Big).
\end{eqnarray*}
By means of independence %and Lemma \ref{l2.4} 
we obtain
\begin{eqnarray*}
\nu_{p-1,0}'(\bar f) &=&
\beta^2u_N^2\sum_{\hat J\in \asp}
(p-1)\big[V_1-4V_2+3V_3\big]+\mathcal{O}\Big(\frac{1}{N}\Big),
\end{eqnarray*}
with

\begin{eqnarray*}
V_1&=& 
\nu_{p-1,0} (\eta_{\hat J}^{1}\eta_{\hat J}^{2}-q^{p-1})
\Bigg[ \nu_{p-1,0} 
\Bigg(\prod_{j\le p-2} \ep_j^1 \ep_j^2\Bigg)-q^{p-1}
 \nu_{p-1,0}\big( \ep_{p-1}^{1}\ep_{p-1}^{2}\big)\Bigg],\\
V_2&=&
\nu_{p-1,0} (\eta_{\hat J}^{1}\eta_{\hat J}^{3}
-q^{p-1})
\Bigg[ \nu_{p-1,0} 
\Bigg(\Big(\prod_{j\le p-2} \ep_j^1 \ep_j^2\ \Big) 
\ep_{p-1}^2 \ep_{p-1}^3\Bigg)\\
&&\quad-q^{p-1}
 \nu_{p-1,0}\big( \ep_{p-1}^{1}\ep_{p-1}^{3}\big)\Bigg],\\
V_3&=&
\nu_{p-1,0} (\eta_{\hat J}^{3}\eta_{\hat J}^{4}-q^{p-1})
\Bigg[ \nu_{p-1,0} 
\Bigg(\Big(\prod_{j\le p-1} \ep_j^1 \ep_j^2\ 
\Big)\ep_{p-1}^3 \ep_{p-1}^4\Bigg)\\
&&\quad-q^{p-1}
 \nu_{p-1,0}\big( \ep_{p-1}^{3}\ep_{p-1}^{4}\big)\Bigg].\\
\end{eqnarray*}
Now, clearly, 
\begin{equation}\label{ecota1}
\nu_{p-1,0}\big( \ep_{p-1}^{1}\ep_{p-1}^{2}\big)
-4 \nu_{p-1,0}\big( \ep_{p-1}^{1}\ep_{p-1}^{3}\big)
+3 \nu_{p-1,0}\big( \ep_{p-1}^{3}\ep_{p-1}^{4}\big)
=0.  
\end{equation}
So, using Lemma \ref{l2.5} together with 
(\ref{ecota1}) we have 
\begin{eqnarray}
\nu_{p-1,0}'(\bar f) &=&
\beta^2u_N^2\sum_{\hat J\in \asp}
(p-1)\nu_{p-1,0} (\eta_{\hat J}^{1}\eta_{\hat J}^{2}-q^{p-1})\nonumber\\
&&\quad\times \Big[q^{p-2}(1-4q+3\hat q_4)\Big]
+\mathcal{O}\Big(\frac{1}{N}\Big).\nonumber
\end{eqnarray}
On the other hand, Lemma \ref{cardasp} implies
$$
u_N^2|\asp|=\frac{p}{2}+\mathcal{O}\Big(\frac{1}{N}\Big).
$$
Then, Corollary \ref{noudos} gives us
\begin{equation}\label{e3.3.11}
\nu_{p-1,0}'(\bar f) =
\frac{\beta^2 p(p-1)}{2}
q^{p-2}(1-4q+3\hat q_4)
\nu_{p-1,0} (R_{1,2}^{p-1}-q^{p-1})
+\mathcal{O}\Big(\frac{1}{N}\Big).
\end{equation}
Invoking inequalities (\ref{e3.3.7'}), (\ref{e3.3.10}) and (\ref{e3.3.11}), 
we now get that
\begin{eqnarray}\label{e3.3.12}
&&\Bigg|\nu (R_{1,2}^{p-1}-q^{p-1})    
- 
\frac{\beta^2 p(p-1)}{2}\nonumber\\
&&\qquad\quad\times q^{p-2}\big[1-4q+3\hat q_4\big]
\nu_{p-1,0} (R_{1,2}^{p-1}-q^{p-1})
\Bigg|\le \frac{K}{N}.
\end{eqnarray}
Finally, using (\ref{e3.5}) and Corollary \ref{c3.5} for $l=1$, we have
\begin{equation}\label{e3.3.13}
\big|\nu(R_{1,2}^{p-1}-q^{p-1})-
\nu_{p-1,0}(R_{1,2}^{p-1}-q^{p-1})\big|\le \frac{K}{N}.
\end{equation}
Then, (\ref{e3.3.12}) and (\ref{e3.3.13}) ensure that for $\beta\le \beta_p$,
(\ref{e3.3.7}) is satisfied.
\hfill
$\Box$

%\vspace{0.1cm}

\section{Study of  the free energy}
Set 
$$p_N(\beta,h,p)=\frac{1}{N}
\be \big[\log Z_N(\beta,h)\big]
=\frac{1}{N} \be \Big[\log
\Big(\sum_{\si\in\ssn}\exp\lp- H_{N,\beta,h}(\si)  \rp
\Big)\Big].$$
This quantity is 
the {\it
expected density of the logarithm of the partition
function}, and sometimes, we will write  $p_N$ instead
of $p_N(\beta,h,p)$. 
The quantity $p_N$ is closely related to the 
{\it free energy} considered by physicists, up to a scaling factor, and we will call it the free energy of our system by usual analogy.

The main aim of this section is to prove the following result.

%\vspace{0.1cm}

\begin{theorem}\label{t4.1}
If $\beta$ satisfies condition {\rm (H)}, we have
%There exists $\beta_p>0$ such that, for any $\beta\le \beta_p$, we have
\begin{eqnarray*}
\lim_{N\uparrow \infty} p_N(\beta,h,p)&=&
\frac{\beta^2}{4}\big[1-pq^{p-1}+(p-1)q^p\big]\\
&&\quad+ \log 2 + \be\Bigg[\log \cosh
\Big[\beta \Big(\frac{p}{2}\Big)^{\frac{1}{2}}
q^{\frac{p-1}{2}}Y+h\Big]\Bigg],
\end{eqnarray*}
where $Y$ is a standard Gaussian random variable and $q$ is the unique solution
to the equation (\ref{qq}).
\end{theorem}
%\vspace{0.1cm}
Before proving Theorem \ref{t4.1},
we will need to introduce some notation and to prove some preliminary results:
consider the function $F:\R_+\times\R_+\times\ou\times \N^*$ defined by
\begin{eqnarray*}
F(\beta,h,q,p)&=&
\frac{\beta^2}{4}\big[1-pq^{p-1}+(p-1)q^p\big]\nonumber\\
&&\quad+\log 2 + \be\Bigg[\log \cosh
\Big[\beta \Big(\frac{p}{2}\Big)^{\frac{1}{2}}
q^{\frac{p-1}{2}}Y+h\Big]\Bigg].
\end{eqnarray*}
We set
\begin{displaymath}
\Phi(\beta,h,p)=F(\beta,h,q,p),
\end{displaymath}
where $q$ satisfies (\ref{qq}). 

%\vspace{0.1cm}

\begin{lemma}\label{l4.2}We have the following two facts
\begin{equation}
{\it  q\ is\ solution\ of\ (\ref{qq})}\ \Rightarrow
\frac{\partial F}{\partial q}(\beta,h,q,p)=0,\label{e4.3}
\end{equation}
\begin{equation}
\frac{\partial \Phi}{\partial \beta}(\beta,h,p)=\frac{\beta}{2}(1-q^p).
\label{e4.4}  \end{equation}
\end{lemma}

%\vspace{0.2cm}

\noindent {\bf Proof:} We first prove  (\ref{e4.3}).
Using integration by parts formula and (\ref{qq}) we obtain
\begin{eqnarray*}
\frac{\partial F}{\partial q}(\beta,h,q,p)&=&
\frac{\beta^2}{4}\big[-p(p-1)q^{p-2}+p(p-1)q^{p-1}\big]\nonumber\\
&&\quad+\beta \frac{p-1}{2}  
\Big(\frac{p}{2}\Big)^{\frac{1}{2}} 
q^{\frac{p-3}{2}}\ 
\be\Bigg[Y \tanh
\Big[\beta \Big(\frac{p}{2}\Big)^{\frac{1}{2}}
q^{\frac{p-1}{2}}Y+h\Big]\Bigg]\\
&=&
\frac{\beta^2}{4}\big[-p(p-1)q^{p-2}+p(p-1)q^{p-1}\big]\nonumber\\
&&\quad+\beta^2
\frac{p(p-1)}{4}
q^{p-2} 
\be\Bigg[ 
\cosh^{-2}
\Big[\beta \Big(\frac{p}{2}\Big)^{\frac{1}{2}}
q^{\frac{p-1}{2}}Y+h\Big]\Bigg]\\
&=&0,
\end{eqnarray*}
which proves (\ref{e4.3}). 
We now show (\ref{e4.4}). This previous result together with
integration by parts and (\ref{qq}) yield

\begin{eqnarray*}
\frac{\partial \Phi}{\partial \beta}(\beta,h,p)&=&
\frac{\partial F}{\partial \beta}\big(\beta,h,q(\beta,h,p),p\big)
+
\frac{\partial F}{\partial q}\big(\beta,h,q(\beta,h,p),p\big)
\ \frac{\partial q}{\partial \beta}(\beta,h,p)\\
&=&
\frac{\beta}{2}\big[1-pq^{p-1}+(p-1)q^{p}\big]\\
&&\quad+\Big(\frac{p}{2}\Big)^{\frac{1}{2}}
q^{\frac{p-1}{2}}\ 
\be\Bigg[Y \tanh
\Big[\beta \Big(\frac{p}{2}\Big)^{\frac{1}{2}}
q^{\frac{p-1}{2}}Y+h\Big]\Bigg]\\
&=&
\frac{\beta}{2}\big[1-pq^{p-1}+(p-1)q^{p}\big]\\
&&\quad+\beta \frac{p}{2}
q^{p-1}\ 
\be\Bigg[ \cosh^{-2}
\Big[\beta \Big(\frac{p}{2}\Big)^{\frac{1}{2}}
q^{\frac{p-1}{2}}Y+h\Big]\Bigg]\\
&=& \frac{\beta}{2} (1-q^p).
\end{eqnarray*}
\hfill
$\Box$

%\vspace{0.2cm}

The following result, that relates the free energy and the overlap, has been proved by Talagrand \cite[Proposition 2.9]{Tp}.
\begin{lemma}\label{l4.4} 
We have
%\begin{equation}\label{e4.7}
$$
\Bigg|\be\Big[\frac{1}{N}\frac{\partial \log Z_N}{\partial \beta}\Big]
-\frac{\beta}{2}\big[1-\be\langle R_{1,2}^p\rangle\big]\Bigg|\le \frac{K}{N}.
$$
%\end{equation}
\end{lemma}

%\vspace{0.2cm}

\noindent
Now we are going to prove the following theorem which 
implies Theorem \ref{t4.1}.

%\vspace{0.1cm}
\begin{theorem}\label{t4.5}
Whenever $\beta$ satisfies {\rm (H)}, we have
%There exists $\beta_p$ such that, for any $\beta \le \beta_p$, we have
%\begin{equation}\label{e4.8}
$$
\big|p_N(\beta,h,p)- \Phi(\beta,h,p)\big|\le \frac{K}{N}.
$$
%\end{equation}
\end{theorem}

%\vspace{0.1cm}

\noindent {\bf Proof:}
We only need to prove that  $p_N(0,h,p)=\Phi(0,h,p)$ and
\begin{equation}\label{e4.9}
\Bigg|\frac{\partial p_N}{\partial \beta}(\beta,h,p)- 
\frac{\partial \Phi}{\partial \beta}(\beta,h,p)\Bigg|\le\frac{K}{N},
\end{equation}
for any $\beta\le \beta_p$. 
For the case $\beta=0$, it is obvious since 
$p_N(0,h,p)=\Phi(0,h,p)=\log (2\cosh h)$.
On the other hand, Lemmas \ref{l4.2},  \ref{l4.4} and
Corollary \ref{cmajoration} imply (\ref{e4.9}).
\hfill
$\Box$

%\vspace{0.1cm}

\section{Almeida-Thouless Theorem}

In this section we prove a  result given in \cite{Al} for the 
Sherrington-Kirkpatrick model. Since the quantity $\Delta_{p-1}^2$ defined above is almost surely positive, it gives a straightforward condition on $\beta$ for the self averaginng to hold, namely that
\begin{equation}\label{defat}
1-\frac{p(p-1)}{2}
q^{p-2}(1-2q+ \hat q_4) \beta^2 >0.
\end{equation}
In fact, this inequality should give the physical limit of the high  temperature region.

\begin{proposition}\label{p5.1}
If $\beta$ satisfies {\rm (H)}, we have
\begin{displaymath}
\Bigg|\nu\Big(
\Delta_{p-1}^2\Big)
-\frac{4(p-1)^2q^{2(p-2)}(1-2q+\hat q_4)}{N\Big(1-\frac{p(p-1)}{2}
q^{p-2}(1-2q+ \hat q_4) \beta^2\Big)}\Bigg|\le \frac{K}{N^{\frac{3}{2}}},
\end{displaymath}
where 
$$\Delta_{p-1}=R_{1,3}^{p-1}-R_{1,4}^{p-1}-R_{2,3}^{p-1}+R_{2,4}^{p-1}.$$
\end{proposition}

\begin{remark}\label{limbeat}
Denote by $\beta_{\mbox{\tiny at}}$ the limit of the region defined by (\ref{defat}), that is
$$
\beta_{\mbox{\tiny at}}=
\beta_{\mbox{\tiny at}}(p)=
\frac{2}{p(p-1)q^{p-2}(1-2q+ \hat q_4)}
=\frac{2}{p(p-1)q^{p-2} E[\cosh^{-4}(Z)]},
$$
where 
$$
Z=\beta \Big(\frac{p}{2}\Big)^{\frac{1}{2}}
q^{\frac{p-1}{2}}Y+h. 
$$
Then, since $q$ tends to $\tanh^2(h)$ as $p\to\infty$, the exponential decay of $q^{p-2}$ implies that
$$
\lim_{p\to\infty}\beta_{\mbox{\tiny at}}(p)=\infty.
$$

\end{remark}

%\vspace{0.1cm}

\noindent {\bf Proof of Proposition \ref{p5.1}:} 
Invoking (\ref{e3.1}) it follows that
\begin{displaymath}
\nu\Big(\Delta_{p-1}^2\Big)
=\frac{2}{p}\nu\Bigg(\Bigg[u_N^2\sum_{J\in Q_{N,1}^p}(\eta_J^1-\eta_J^2)
(\eta_J^3-\eta_J^4)+V_N\Bigg]\Delta_{p-1}\Bigg),
\end{displaymath}
with 
\begin{equation}\label{e5.1}
|V_N|\le \frac{K}{N}.
\end{equation}
The symmetry property  implies now that
\begin{displaymath}
\nu\Big(\Delta_{p-1}^2\Big)
=\frac{2}{p}\nu\Big(\big[u_N^2 |Q_{N,1}^p|(\ep^1-\ep^2)
(\ep^3-\ep^4)+V_N\big]\Delta_{p-1}\Big),
\end{displaymath}
where we have used the notation $\ep^l=\ep^l_1\times\cdots \times \ep_{p-1}^l$
for $l=1,2,3,4$. 

First of all we will check that $V_N$ gives raise to a negligible term: indeed,
Corollary  \ref{c3.5} yields
\begin{equation}\label{e*star}
\nu\Big(\Delta_{p-1}^2 \Big)
\le \frac{K}{N}.
\end{equation}
Then, estimates (\ref{e5.1}), (\ref{e*star}) together with Lemma
\ref{cardbn} give 
$$\nu\Big(\Delta_{p-1}^2\Big)
=\nu\Big((\ep^1-\ep^2)
(\ep^3-\ep^4)
\Delta_{p-1}\Big)+ \mathcal{O}\Big(\frac{1}{N^{\frac{3}{2}}}\Big).
$$

For $(l,l')\in \{(1,3), (1,4), (2,3), (2,4)\}$,  we will now decompose $R_{l,l^{'}}$ into a part involving only the $N-(p-1)$ first spins on one hand, and a remaining term on the other hand, as follows:
\begin{displaymath}
R_{l,l'}^{p-1}=
\Bigg[R_{l,l'}^*+\frac{1}{N}\sum_{j=1}^{p-1}\ep_j^l \ep_j^{l'}\Bigg]^{p-1}
=\sum_{k=0}^{p-1}{p-1 \choose k} \big(R_{l,l'}^*\big)^k 
\Bigg[\frac{1}{N}\sum_{j=1}^{p-1}\ep_j^l \ep_j^{l'}\Bigg]^{p-1-k},
\end{displaymath}
where
$$R_{l,l'}^*=\frac{1}{N} \sum_{j=1}^{N-(p-1)} \sigma_j^l \sigma_j^{l'}.$$
Now, set
$$\Delta_{p-1}^{*}= \big(R_{1,3}^*\big)^{p-1}
-\big(R_{1,4}^*\big)^{p-1}-    \big(R_{2,3}^*\big)^{p-1}    
+ \big(R_{2,4}^*\big)^{p-1}.$$    
Using this decomposition we have
\begin{equation}\label{e5.2}
\nu\Big(\Delta_{p-1}^2\Big)
=\nu\Bigg((\ep^1-\ep^2)
(\ep^3-\ep^4)
\Bigg[\Delta_{p-1}^{*}+\sum_{k=0}^{p-2} {p-1 \choose k} y_k \Bigg]\Bigg)
+\mathcal{O}\Big(\frac{1}{N^{\frac{3}{2}}}\Big)  ,
 \end{equation}
with
\begin{eqnarray*}
y_k&=&
\big(R_{1,3}^*\big)^k 
\Bigg[\frac{1}{N}\sum_{j=1}^{p-1}\ep_j^1 \ep_j^{3}\Bigg]^{p-1-k}
-\big(R_{1,4}^*\big)^k 
\Bigg[\frac{1}{N}\sum_{j=1}^{p-1}\ep_j^1 \ep_j^{4}\Bigg]^{p-1-k}\\
&&- \big(R_{2,3}^*\big)^k 
\Bigg[\frac{1}{N}\sum_{j=1}^{p-1}\ep_j^2 \ep_j^{3}\Bigg]^{p-1-k}
+ \big(R_{2,4}^*\big)^k 
\Bigg[\frac{1}{N}\sum_{j=1}^{p-1}\ep_j^2 \ep_j^{4}\Bigg]^{p-1-k}.
\end{eqnarray*}
We will deal now with the different terms appearing in (\ref{e5.2}) separately, and we start with all the terms containing $y_k$. By the construction of $y_k$ it
is easily checked that
for any $k \in
\{0,\dots,p-2\}$,  
$$|y_k | \le \frac{K}{N^{p-1-k}}.$$
These bounds and inequality (\ref{e3.5}) in Corollary \ref{c3.6} 
yield, for any $k \in
\{0,\dots,$ $p-2\}$,  
\begin{equation}\label{e5.3}\Big|\nu 
\Big((\ep^1-\ep^2)
(\ep^3-\ep^4) y_k\Big)
- \nu _{p-1,0}\Big((\ep^1-\ep^2)
(\ep^3-\ep^4) y_k\Big)\Big| \le \frac{K}{N^{\frac{3}{2}}}.
\end{equation}
Hence, we are reduced to study terms of the form $\nu _{p-1,0}((\ep^1-\ep^2)
(\ep^3-\ep^4) y_k)$. 
Since for any $k\in \{0,\dots,p-3\}$,
$$
\Big| \nu _{p-1,0}\big((\ep^1-\ep^2)
(\ep^3-\ep^4) y_k\big)\Big|\le \frac{K}{N^2},
$$
%\end{equation}
we only have to deal with the term involving $y_{p-2}$. But the definition of $y_{p-2}$
implies that
\begin{equation}\label{e5.6}\begin{array}{l} 
\displaystyle\nu _{p-1,0}\Big((\ep^1-\ep^2)
(\ep^3-\ep^4) y_{p-2}\Big)=\frac{1}{N}\sum_{j=1}^{p-1}
\nu _{p-1,0}\Big((\ep^1-\ep^2)
(\ep^3-\ep^4)\\  
\quad\times\Big[\big(R_{1,3}^*\big)^{p-2}\ep_j^1 \ep_j^3
-\big(R_{1,4}^*\big)^{p-2}\ep_j^1 \ep_j^4
-\big(R_{2,3}^*\big)^{p-2}\ep_j^2 \ep_j^3
+\big(R_{2,4}^*\big)^{p-2}\ep_j^2 \ep_j^4\Big]\Big). 
\end{array}\end{equation}
Observe now the first term of the right hand side of (\ref{e5.6}).  Using the independence, 
the meaning of $\ep^l$
and Lemma \ref{l2.5}, we obtain
\begin{equation}\label{e5.7}
\begin{array}{l}
\displaystyle\nu _{p-1,0}\Big((\ep^1-\ep^2)
(\ep^3-\ep^4)\big(R_{1,3}^*\big)^{p-2}\ep_j^1 \ep_j^3\Big)\\ 
\displaystyle \quad =\nu _{p-1,0}\Big(\big(R_{1,3}^*\big)^{p-2}\Big)
\nu _{p-1,0}\Big((\ep^1-\ep^2)
(\ep^3-\ep^4)\ep_{p-1}^1 \ep_{p-1}^3\Big)\\ 
\displaystyle \quad = \nu _{p-1,0}\Big(\big(R_{1,3}^*\big)^{p-2}\Big)\\
\displaystyle \qquad \times 
\nu _{p-1,0}\Bigg(\prod_{j=1}^{p-2}\ep_j^1 \ep_j^3
-2 \Bigg[\prod_{j=1}^{p-2}\ep_j^1 \ep_j^4\Bigg] \ep_{p-1}^3 \ep_{p-1}^4 
+\ep^2 \ep^4 \ep_{p-1}^1 \ep_{p-1}^3\Bigg)\\ 
\displaystyle \quad = \nu _{p-1,0}\Big(\big(R_{1,3}^*\big)^{p-2}\Big)  
\ q^{p-2} (1-2q+\hat q_4) + \mathcal{O}\Big(\frac{1}{N}\Big).
\end{array}\end{equation}
Let us study now $\nu _{p-1,0}\Big( \big(R_{1,3}^*\big)^{p-2}\Big)$.
The inequality (\ref{e3.5}) in Corollary \ref{c3.6} yields
\begin{equation}\label{e5.8}
\Big|\nu\Big(\big(R_{1,3}^*\big)^{p-2}\Big)
-\nu _{p-1,0}\Big(\big(R_{1,3}^*\big)^{p-2}\Big)\Big|\le 
\frac{K}{N^\frac12},
\end{equation}
the definition of 
$R_{1,3}^*$ implies
\begin{equation}\label{e5.9}
\Big|\big(R_{1,3}^*\big)^{p-2}
-R_{1,3}^{p-2}\Big|\le \frac{K}{N},
\end{equation}
and on the other hand, Corollary \ref{cmajoration} gives us
\begin{equation}\label{e5.10}
\Big|\nu\Big(R_{1,3}^{p-2}-q^{p-2}\Big)\Big|\le \frac{K}{N}.
\end{equation}
So, putting togheter (\ref{e5.8})-(\ref{e5.10}) we have
\begin{equation}\label{XX}
\nu _{p-1,0}\Big( \big(R_{1,3}^*\big)^{p-2}\Big)=q^{p-2}+\mathcal{O}\Big(\frac{1}{N^{\frac{1}{2}}}\Big).
\end{equation}
The other terms in (\ref{e5.6}) can be studied in a similar way. Then, putting together 
(\ref{e5.3}), (\ref{e5.6}), (\ref{e5.7}) and  (\ref{XX}) we easily get
\begin{multline}\label{e5.11}
\nu\Bigg((\ep^1-\ep^2)
(\ep^3-\ep^4)\sum_{k=0}^{p-2} {p-1 \choose k} y_k \Bigg)\\
= \frac{4(p-1)^2q^{2(p-2)}(1-2q+\hat q_4)}{N}
+ \mathcal{O}\Big(\frac{1}{N^{\frac{3}{2}}}\Big).
\end{multline}

Now we will deal with the term containing $\Delta_{p-1}^{*}$ in 
(\ref{e5.2}), that means with
$$
\nu\Big((\ep^1-\ep^2)
(\ep^3-\ep^4) \Delta_{p-1}^{*} \Big).
$$
To that purpose, we will use 
(\ref{e3.6'}) in Corollary  \ref{c3.6}. 
First of all, note that
\begin{equation}\label{noupa}
\vert \Delta_{p-1} - \Delta_{p-1}^* \vert \le \frac{K}{N}.
\end{equation}
Then estimate (\ref{e*star}) yields 
\begin{equation}\label{ee*star}
\nu\Big(\big(\Delta_{p-1}^*\big)^2 \Big)
\le \frac{K}{N}.
\end{equation}
Since the  independence ensures
\begin{equation}\label{e5.12}
\displaystyle\nu_{p-1,0}\Big((\ep^1-\ep^2)
(\ep^3-\ep^4) \Delta_{p-1}^{*} \Big)
= \nu_{p-1,0}\big(
(\ep^1-\ep^2)(\ep^3-\ep^4)\big) \nu_{p-1,0}\big(\Delta_{p-1}^{*} \big)=0,
\end{equation}
we have,  using (\ref{ee*star}) 
\begin{equation}\label{e5.14}%\begin{array}{l}
\displaystyle\nu \Big((\ep^1-\ep^2)
(\ep^3-\ep^4) \Delta_{p-1}^{*} \Big)%\\
%\displaystyle\quad\quad\quad
= \nu_{p-1,0}'\big(
(\ep^1-\ep^2)(\ep^3-\ep^4) \Delta_{p-1}^{*} \big)
+ \mathcal{O}\Big(\frac{1}{N^{\frac{3}{2}}}\Big).
%\end{array}
\end{equation}
Moreover, Proposition \ref{pderivnu} yields
\begin{equation}\label{e5.15}
\nu_{p-1,0}'\Big((\ep^1-\ep^2)
(\ep^3-\ep^4) \Delta_{p-1}^{*} \Big)=\beta^2u_N^2\sum_{J\in Q_{N,p-1}^p}
\big(D^J_1+D^J_2+D^J_3\big),
\end{equation}
with
\begin{eqnarray*}
D_1^J&=& \nu_{p-1,0}\Big( \Delta_{p-1}^{*}(\ep^1-\ep^2)(\ep^3-\ep^4)
\sum_{1\le l<l'\le 4}(\eta_J^l\eta_J^{l'}-q^{p-1}) 
\ep_J^l\ep_J^{l'}\Big),\\
D_2^J&=&-4 \nu_{p-1,0}\Big( \Delta_{p-1}^{*} (\ep^1-\ep^2)(\ep^3-\ep^4)
\sum_{l\le 4}(\eta_J^l\eta_J^{5}-q^{p-1}) 
\ep_J^l\ep_J^{5}\Big),\\
D_3^J&=&10 \nu_{p-1,0}\Big( \Delta_{p-1}^{*}(\ep^1-\ep^2)(\ep^3-\ep^4)
(\eta_J^5\eta_J^{6}-q^{p-1}) 
\ep_J^5\ep_J^{6}\Big).
\end{eqnarray*}
We  can check that 
$$D_2^J\equiv 0,\qquad {\rm for\ any}\  J \in Q_{N,p-1}^{p-1}.$$
Indeed, for instance, when $l=1$, we get that
\begin{eqnarray*} 
& &\nu_{p-1,0}\Big( \Delta_{p-1}^{*}
(\ep^1-\ep^2)(\ep^3-\ep^4)(\eta_J^1\eta_J^{5}-q^{p-1}) 
\ep_J^1\ep_J^{5}\Big)\\
& &\qquad =
\nu_{p-1,0}\Big( \Delta_{p-1}^{*} (\eta_J^1\eta_J^{5}-q^{p-1})\Big)
\nu_{p-1,0}\Big((\ep^1-\ep^2)(\ep^3-\ep^4) 
\ep_J^1\ep_J^{5}\Big).
\end{eqnarray*}
Moreover,
\begin{multline*}
\nu_{p-1,0}\Big((\ep^1-\ep^2)(\ep^3-\ep^4) \ep_J^1\ep_J^{5}\Big)\\
=\nu_{p-1,0}\Big((\ep^1-\ep^2)\ep^3\ep_J^1\ep_J^{5}\Big)
-\nu_{p-1,0}\Big((\ep^1-\ep^2)\ep^4\ep_J^1\ep_J^{5}\Big)
=0.
\end{multline*}
This kind of argument, that will be repeated all along the remainder of the paper, will be referred to as {\it symmetry among the different copies of $G_N$}. Now the cases $l=2,3,4$ in $D_2^J$ can be studied with the same method, and furthermore, by similar arguments,
$$D_3^J\equiv 0,\qquad {\rm for\ any}\  J \in Q_{N,p-1}^{p-1}.$$
Thus, it only remains to deal with $D_1^J$.
The summatory of $D_1^J$ contains the couples
$$(l,l')\in \{(1,2), (1,3), (1,4), (2,3), (2,4), (3,4)\}.$$
We start studying the couple $(1,2)$. 
Along the same lines as before, using independence and symmetry between copies of $G_N$, we easily get that
$$\begin{array}{l}
\displaystyle\nu_{p-1,0}\Big( \Delta_{p-1}^{*}
(\ep^1-\ep^2)(\ep^3-\ep^4)(\eta_J^1\eta_J^{2}-q^{p-1}) 
\ep_J^1\ep_J^{2}\Big)\\[5mm]
\displaystyle\quad\quad =
\nu_{p-1,0}\Big( \Delta_{p-1}^{*}(\eta_J^1\eta_J^{2}-q^{p-1})\Big)
\nu_{p-1,0}\Big((\ep^1-\ep^2)(\ep^3-\ep^4) 
\ep_J^1\ep_J^{2}\Big)\\[5mm]
\quad\quad=0,
\end{array}$$
and the same argument can be applied to $(3,4)$.
Consider now the couple $(1,3)$ and the summatory of 
$\sum_{J\in Q_{N,p-1}^p}D_1^J$. 
Lemma \ref{darrer}, 
(\ref{ee*star}), the properties of independence 
 and symmetry among all the $\ep_j$ imply
$$\begin{array}{l}
\displaystyle u_N^2 
\sum_{J\in Q_{N,p-1}^p} \nu_{p-1,0}\Big( \Delta_{p-1}^{*}
(\ep^1-\ep^2)(\ep^3-\ep^4)(\eta_J^1\eta_J^{3}-q^{p-1}) 
\ep_J^1\ep_J^{3}\Big),\\[5mm]
%\displaystyle 
%\qquad=u_N^2 \sum_{J\in Q_{N,p-1}}
%\nu_{p-1,0}\Big( \Delta_{p-1}^{*}(\eta_J^1\eta_J^{3}-q^{p-1})\Big)%\\[5mm]
%\nu_{p-1,0}\Big((\ep^1-\ep^2)(\ep^3-\ep^4) 
%\ep_J^1\ep_J^{3}\Big),\\[5mm]
\displaystyle 
\qquad= u_N^2 \sum_{\hat J\in \asp} \sum_{l=1}^{p-1}
\nu_{p-1,0}\Big( \Delta_{p-1}^{*}(\eta_{\hat J}^1\eta_{\hat J}^{3}-q^{p-1})\Big)\\[5mm]
\qquad\qquad\qquad\times\displaystyle
\nu_{p-1,0}\Big((\ep^1-\ep^2)(\ep^3-\ep^4) 
\ep_l^1\ep_l^{3}\Big)+\mathcal{O}\Big(\frac{1}{N^{\frac{3}{2}}}\Big),\\[5mm]
\displaystyle 
\qquad= u_N^2 \sum_{\hat J\in \asp} (p-1)
\nu_{p-1,0}\Big( \Delta_{p-1}^{*}(\eta_{\hat J}^1\eta_{\hat J}^{3}-q^{p-1})\Big)\\[5mm]
\qquad\qquad\qquad\times\displaystyle
\nu_{p-1,0}\Big((\ep^1-\ep^2)(\ep^3-\ep^4) 
\ep_{p-1}^1\ep_{p-1}^{3}\Big)+\mathcal{O}\Big(\frac{1}{N^{\frac{3}{2}}}\Big).
\end{array}$$
Notice also that by means of  Lemma \ref{l2.5}, we have
$$ 
\nu_{p-1,0}\Big((\ep^1-\ep^2)(\ep^3-\ep^4) 
\ep_{p-1}^1\ep_{p-1}^{3}\Big)
= q^{p-2}(1-2q+\hat q_4)+ 
\mathcal{O}\Big(\frac{1}{N}\Big).
$$
Let us handle now the remaining couples in $D_1^J$: it is easily seen that
the couple $(2,4)$ has the same structure as $(1,3)$ and the couples
$(2,3)$ and $(1,4)$ also have the same structure as $(1,3)$
but with the opposite sign.
Plugging  now (\ref{e5.11}), (\ref{e5.12}), (\ref{e5.14}), and
 the study of $D_1^J, D_2^J$, 
$D_3^J$ in (\ref{e5.15}) into (\ref{e5.2}), 
we end up with
\begin{eqnarray*}
\nu\Big(
\Delta_{p-1}^2\Big)&&=
\frac{4(p-1)^2q^{2(p-2)}(1-2q+\hat q_4)}{N}
+(p-1)q^{p-2}(1-2q+\hat q_4) \\
&&  \times\beta^2 u_N^2 \sum_{\hat J\in \asp}
\nu_{p-1,0}\Big( \Delta_{p-1}^{*}(\eta_{\hat J}^1\eta_{\hat J}^{3}
-\eta_{\hat J}^1\eta_{\hat J}^{4}-
\eta_{\hat J}^2\eta_{\hat J}^{3}
+\eta_{\hat J}^2\eta_{\hat J}^{4})\Big)\\
&&+\mathcal{O}\Big(\frac{1}{N^{\frac{3}{2}}}\Big).
\end{eqnarray*}
Moreover, (\ref{noupa}),  Corollary \ref{noudos} and estimate (\ref{e*star}) imply that
\begin{eqnarray*}
\nu\Big(
\Delta_{p-1}^2\Big)&=&
\frac{4(p-1)^2q^{2(p-2)}(1-2q+\hat q_4)}{N}
+\beta^2\frac{p}{2}(p-1)q^{p-2}\\
&& \ \times (1-2q+\hat q_4) 
\nu_{p-1,0}\Big( \Delta_{p-1}^2
\Big)+\mathcal{O}\Big(\frac{1}{N^{\frac{3}{2}}}\Big).
\end{eqnarray*}
Finally, inequality (\ref{e3.5}) in Corollary \ref{c3.6} and an
upper bound for $\Delta_{p-1}^4$ obtained from Corollary
\ref{c3.5}, and similar to the one obtained at (\ref{e*star}), ensure
\begin{eqnarray*}
\nu\Big(
\Delta_{p-1}^2\Big)&=&
\frac{4(p-1)^2q^{2(p-2)}(1-2q+\hat q_4)}{N}
+\beta^2\frac{p}{2}(p-1)q^{p-2}\\
&& \ \times (1-2q+\hat q_4) 
\nu\Big( \Delta_{p-1}^2
\Big)+\mathcal{O}\Big(\frac{1}{N^{\frac{3}{2}}}\Big).
\end{eqnarray*}
This last estimate immediately gives the desired result.
\hfill
$\Box$

%\vspace{0.1cm}

\section{Second Moment Computations}

The results of this section will be crucial to obtain our
Central Limit Theorems: in fact, the relations we will obtain in the next Section for any power $k$ will be a mere elaboration of the ones obtained here for $k=2$.

For any $\J=(i_1,\ldots,i_{p-1}) \in N_{p-1}$ (see Definition \ref{setnr})  we define $\eta_{\J}$ as
\begin{equation}\label{e6.3'}
\eta_{\J}=\prod_{j=1}^{p-1} \sigma_{i_j}.
\end{equation}
Then set
$$T(\eta)=\{\eta_{\J};\, \J \in N_{p-1}\}$$ 
and 
$$
b=\langle T(\eta) \rangle=\lcl
\langle T(\eta_{\J}) \rangle ;\,\J \in N_{p-1}\rcl.
$$
Let us define now
\begin{equation}\begin{array}{ll}
T_{l,l'}&=\displaystyle\frac{1}{N^{p-1}}\big(T(\eta^l)-b\big)\cdot
\big(T(\eta^{l'})-b\big)
=\frac{1}{N^{p-1}}\sum_{\J \in N_{p-1}}(\eta^l_{\J}-b_{\J})
(\eta^{l'}_{\J}-b_{\J})\\%\label{e6.1}\\
T_l&=\displaystyle\frac{1}{N^{p-1}}\big(T(\eta^l)-b\big)\cdot b
=\frac{1}{N^{p-1}}\sum_{\J \in N_{p-1}}
(\eta^l_{\J}-b_{\J}) b_{\J}\\%\label{e6.2}\\
T&=\displaystyle\frac{b\cdot b}{N^{p-1}}-q^{p-1}.%\label{e6.3}
\end{array}\nonumber\end{equation}
Equality (\ref{etll})  allows us to reconstruct the quantities $R_{l,l'}^{p-1}$
by the formula 
\begin{equation}\label{e6.4}
T_{l,l'}+T_l+T_{l'}+T=R_{l,l'}^{p-1}-q^{p-1}.
\end{equation}
We finish this introduction with some more notation that we will use all along this section. First, set
$$\hat T(\eta^l)=
\{\eta_{\J}^l, \J \in A_N^{p-1}\}.$$
On the other hand, note that the set 
$A_{N}^{p-1}$ can be decomposed into three disjoint sets 
 as follows
\begin{equation}\label{e6.12}
A_{N}^{p-1}=\asp \cup   \tilde Q_{N,p-1}^{p-1} \cup
 \bar Q_{N,p-1}^{p-1}.
 \end{equation}
 
We are now ready to estimate the second moment of $T_{l,l^{'}}$, $T_l$ and $T$.
\begin{proposition}\label{p6.1}
If $\beta$ satisfies {\rm (H)}, we have
\begin{equation}\label{e6.5}
\Big|\nu\big(T_{1,2}^2\big)-\frac{A^2}{N}\Big|\le \frac{K}{N^{\frac{3}{2}}}
\end{equation}
where
\begin{equation} \label{e6.6}
A^2=\frac{(p-1)^2 q^{2(p-1)} (1-2q+\hat q_4)}{1-\frac{p(p-1)}{2} 
q^{p-2} (1-2q+\hat q_4)\beta^2}.
\end{equation}
\end{proposition}

%\vspace{0.2cm}

\noindent {\bf Proof:} 
As in the proof of Lemma \ref{nouu},  we  will decompose $T_{1,2 }$ into two terms 
\begin{equation}\label{e6.7}
T_{1,2 }=S_{1,2}+V_{N}^{1,2},
\end{equation}
with
\begin{eqnarray*}
S_{1,2}&=&
\frac{(p-1)!}{N^{p-1}}
\sum_{\J \in A_N^{p-1}}
(\eta^1_{\J}-b_{\J})(\eta^{2}_{\J}-b_{\J}),\\
V_{N}^{1,2}&=&
\frac{1}{N^{p-1}}
\sum_{{\J}\in \bar N_{p-1}^c}
(\eta^1_{\J}-b_{\J})(\eta^{2}_{\J}-b_{\J}),
\end{eqnarray*}
where $\eta_{\J}$ for ${\J} \in A_{N}^{p-1}$ 
or ${\J} \in \bar N_{p-1}^c$ (see the Appendix) are defined as in (\ref{e6.3'}). 
Then, it suffices to study $\nu(S_{1,2}^2)$. Indeed, from Lemma \ref{cardnp}
we have
\begin{equation}\label{YY}
\vert V_{N}^{1,2} \vert \le \frac{K}{N}.
\end{equation}

On the other hand, using the symmetry property, we get 
\begin{eqnarray}
\nu\left(S_{1,2}^2\right)&=&
\left(\frac{(p-1)!}{N^{p-1}} \right)^2
\nu\Big[\left(\hat T(\eta^1)-\hat T(\eta^2)\right) \cdot
\left(\hat T(\eta^3)-\hat T(\eta^4)\right)\nonumber\\
&& \quad \times
\left(\hat T(\eta^1)-\hat T(\eta^5)\right) \cdot
\left(\hat T(\eta^3)-\hat T(\eta^6)\right)\Big]\nonumber\\
&=&\frac{[(p-1)!]^2}{N^{2(p-1)}}
\big| A_{N}^{p-1}\big| 
\nu\Big[(\ep^1-\ep^2) (\ep^3-\ep^4)
\left(\hat T(\eta^1)-\hat T(\eta^5)\right)\nonumber\\
&& \quad \cdot
\left(\hat T(\eta^3)-\hat T(\eta^6)\right)\Big].\label{e6.s12}
\end{eqnarray}
So, Lemma \ref{cardasp}, Corollary \ref{c3.6}, Lemma \ref{nouu} and Proposition \ref{p5.1} imply
$$
\nu\left(S_{1,2}^2\right)=
\mathcal{O}\left(\frac{1}{N}\right).
$$ 
%Using the same ideas we can also prove that
%\begin{equation}\label{ZZ}
%\nu\left(S_{1,2}^4\right)=
%\mathcal{O}\left(\frac{1}{N^2}\right).
%\end{equation}
Then, we clearly have
\begin{equation}\label{e6.11}
\nu\left(T_{1,2}^2\right)=\nu\left(S_{1,2}^2\right)+
\mathcal{O}\left(\frac{1}{N^{\frac{3}{2}}}\right).
\end{equation}

Let us study now $\nu(S_{1,2}^2).$ Relation (\ref{e6.s12}), 
Lemma \ref{cardasp}
and the decomposition of $A_N^{p-1}$ (see (\ref{e6.12})) allow us to write
\begin{equation}\label{e6.18'}
\nu\left(S_{1,2}^2\right)\ =\ M_1\ +\ M_2\ + \ M_3 +
\mathcal{O}\left(\frac{1}{N^{\frac{3}{2}}}\right),
\end{equation} 
with
\begin{eqnarray*}
M_1&=&\frac{(p-1)!}{N^{p-1}} 
\nu\Bigg[(\ep^1-\ep^2) (\ep^3-\ep^4)
\sum_{{\J} \in \asp}
(\eta^1_J-\eta_{\J}^5)(\eta_{\J}^3-\eta^6_{\J})\Bigg],\\
M_2&=&\frac{(p-1)!}{N^{p-1}}
\nu\Bigg[(\ep^1-\ep^2) (\ep^3-\ep^4)
\sum_{{\J}\in \tilde Q_{N,p-1}^{p-1}}
(\eta^1_{\J}-\eta_{\J}^5)(\eta_{\J}^3-\eta^6_{\J})\Bigg],\\
M_3&=&\frac{(p-1)!}{N^{p-1}}
\nu\Bigg[(\ep^1-\ep^2) (\ep^3-\ep^4)
\sum_{{\J}\in \bar Q_{N,p-1}^{p-1}}
(\eta^1_{\J}-\eta_{\J}^5)(\eta_{\J}^3-\eta^6_{\J})\Bigg].\\
\end{eqnarray*}

The term $M_3$ is easily handled: Lemma \ref{cardbn} clearly yields
\begin{equation}\label{e6.19}
\big|M_3\big| \le \frac{K}{N^2}.
\end{equation}

We deal now with $M_2$. By symmetry, we have
\begin{eqnarray*}
M_2&=&\frac{(p-1)!}{N^{p-1}} (p-1) 
\nu\Bigg[(\ep^1-\ep^2) (\ep^3-\ep^4)\\
&&\quad\times \sum_{{\J}\in A_{N-(p-1)}^{p-2}}
(\eta^1_{\J}\ep^1_{p-1}-\eta_{\J}^5\ep^5_{p-1})
(\eta_{\J}^3\ep^3_{p-1}-\eta^6_{\J}\ep^6_{p-1})\Bigg].
\end{eqnarray*}
Corollary \ref{c3.6} and  the independence
ensure that, up to a $N^{-3/2}$ term, we have
\begin{eqnarray*}
M_2&=&\frac{(p-1)!}{N^{p-1}} (p-1)\\ 
&&\quad\times \Bigg[
\nu_{p-1,0}\big((\ep^1-\ep^2) (\ep^3-\ep^4)\ep^1_{p-1} \ep^3_{p-1}\big)
\nu_{p-1,0}\Bigg(\sum_{{\J}\in A_{N-(p-1)}^{p-2}}
\eta^1_{\J} \eta_{\J}^3\Bigg)\\
&&\quad - 2
\nu_{p-1,0}\big((\ep^1-\ep^2) (\ep^3-\ep^4)\ep^1_{p-1} \ep^6_{p-1}\big)
\nu_{p-1,0}\Bigg(\sum_{{\J}\in A_{N-(p-1)}^{p-2}}
\eta^1_{\J} \eta_{\J}^6\Bigg)\\
&&\quad +
\nu_{p-1,0}\big((\ep^1-\ep^2) (\ep^3-\ep^4)\ep^5_{p-1} \ep^6_{p-1}\big) 
\nu_{p-1,0}\Bigg(\sum_{{\J}\in A_{N-(p-1)}^{p-2}}\eta^5_{\J} 
\eta_{\J}^6\Bigg)\Bigg] 
%+\mathcal{O}\left(\frac{1}{N^{\frac{3}{2}}}\right).
\end{eqnarray*}
Let us compute all the terms  of the last equality: 
on one hand, Lemma \ref{l2.5} implies
\begin{displaymath}
\nu_{p-1,0}\big((\ep^1-\ep^2) (\ep^3-\ep^4)\ep^1_{p-1} \ep^3_{p-1}\big) 
=q^{p-2} (1-2q + \hat q_4)+ \mathcal{O}\left(\frac{1}{N}\right),
\end{displaymath}
and an easy symmetry argument yields
\begin{displaymath}
\nu_{p-1,0}\big((\ep^1-\ep^2) (\ep^3-\ep^4)\ep^1_{p-1} \ep^6_{p-1}\big)
= \nu_{p-1,0}\big((\ep^1-\ep^2) (\ep^3-\ep^4)\ep^5_{p-1} \ep^6_{p-1}\big)
=0.           
\end{displaymath}
On the other hand,  an obvious extension of (\ref{equasif}) in Lemma \ref{nouu}
shows that
\begin{displaymath}
R_{1,2}^{p-2}=\frac{(p-2)!}{N^{p-2}}
\sum_{{\J} \in A_{N-(p-1)}^{p-2}}\eta^1_{\J} \eta_{\J}^2
+\mathcal{O}\left(\frac{1}{N}\right).
\end{displaymath}
Then, Corollaries \ref{c3.6} and \ref{cmajoration}
ensure that
\begin{equation}\label{e6.20}
M_2= \frac{1}{N}(p-1)^2 q^{2(p-2)} (1-2q + \hat q_4)     
+\mathcal{O}\left(\frac{1}{N^{\frac{3}{2}}}\right).
\end{equation}
It only remains to study $M_1$. Set
\begin{displaymath}
f=
\frac{(p-1)! }{N^{p-1}}(\ep^1-\ep^2) (\ep^3-\ep^4)
\sum_{{\J}\in \asp}
(\eta^1_{\J}-\eta_{\J}^5)(\eta_{\J}^3-\eta^6_{\J}).
\end{displaymath}
Then
$$M_1= \nu(f).$$
Corollary \ref{noudos} and Proposition \ref{p5.1}  implies
$$\nu^{\frac{1}{2}}\big(f^2\big)\le \frac{K}{N^{\frac{1}{2}}}.$$
Furthermore, since by  symmetry between the different copies of $G_N$, 
\break $\nu_{p-1,0}(f)=0$, Corollary \ref{c3.6} implies
$$\nu(f)= \nu'_{p-1,0}(f) 
+\mathcal{O}\left(\frac{1}{N^{\frac{3}{2}}}\right).$$
Let us now compute $\nu_{p-1,0}'(f)$. Proposition \ref{pderivnu} yields
$$
\nu_{p-1,0}'(f)
=\beta^2u_N^2\sum_{J'\in Q_{N,p-1}^p}
\big(F^{J'}_1+F^{J'}_2+F^{J'}_3\big),
$$
with
\begin{eqnarray*}
F_1^{J'}&=& 
\nu_{p-1,0}\Bigg( f
\sum_{1\le l<l'\le 6}(\eta_{J'}^l\eta_{J'}^{l'}-q^{p-1}) 
\ep_{J'}^l\ep_{J'}^{l'}\Bigg),\\
F_2^{J'}&=&-6 \nu_{p-1,0}\Bigg( f
\sum_{l\le 6}(\eta_{J'}^l\eta_{J'}^{7}-q^{p-1}) 
\ep_{J'}^l\ep_{J'}^{7}\Bigg),\\
F_3^{J'}&=&21 \nu_{p-1,0}\Big( f
(\eta_{J'}^7\eta_{J'}^{8}-q^{p-1}) 
\ep_{J'}^7\ep_{J'}^{8}\Big).
\end{eqnarray*}
Using again the symmetry argument (again among the different copies of $G_N$), the quantities $F^{J'}_2$, $F^{J'}_3$
and all the couples of $F^{J'}_1$ except for
$\{(1,3),$ $(1,4),$  $(2,3),$ $ (2,4)\}$ vanish.
Now, we will analize the couple $(1,3)$ of $F^{J'}_1$.
We follow the ideas given in the proof of Proposition \ref{p5.1}.
Properties of  independence and symmetry, Lemma \ref{darrer}, Proposition \ref{p3.2} and Lemma \ref{l2.5}
ensure

\begin{displaymath}\begin{array}{ll}
&\displaystyle\frac{(p-1)!}{N^{p-1}}  
\beta^2 u_N^2 \sum_{J'\in Q_{N,p-1}^p} \sum_{{\J} \in \asp}  
\nu_{p-1,0}\Big((\ep^1-\ep^2) (\ep^3-\ep^4) \ep_{J'}^1\ep_{J'}^{3}\\[3mm]
&\ \displaystyle\quad \times 
(\eta^1_{\J}-\eta_{\J}^5)(\eta_{\J}^3-\eta^6_{\J})
(\eta_{J'}^1\eta_{J'}^{3}-q^{p-1})\Big)\\[3mm]
&=\displaystyle\frac{(p-1)!}{N^{p-1}} (p-1) 
\beta^2 u_N^2 \\[3mm]
&\ \displaystyle\quad \times 
\sum_{{\J}'\in \asp} \sum_{{\J}\in \asp}  
\nu_{p-1,0}\Big((\ep^1-\ep^2) (\ep^3-\ep^4) \ep_{p-1}^1\ep_{p-1}^{3}\Big)\\[3mm]
&\ \displaystyle\quad \times 
\nu_{p-1,0}\Big((\eta^1_{\J}-\eta_{\J}^5)(\eta_{\J}^3-\eta^6_{\J})
(\eta_{{\J}'}^1\eta_{{\J}'}^{3}-q^{p-1})\Big)
+\mathcal{O}\left(\frac{1}{N^{\frac{3}{2}}}\right)\\[3mm]
&=\displaystyle\frac{(p-1)!}{N^{p-1}} (p-1) 
\beta^2 u_N^2 
q^{p-2} (1-2q + \hat q_4)\\[3mm]
&\ \displaystyle\times\sum_{{\J}'\in \asp} \sum_{{\J}\in \asp}  
\nu_{p-1,0}\Big((\eta^1_{\J}-\eta_{\J}^5)(\eta_{\J}^3-\eta^6_{\J})
(\eta_{{\J}'}^1\eta_{{\J}'}^{3}-q^{p-1})\Big)\\[3mm]
&\quad\displaystyle
+\mathcal{O}\left(\frac{1}{N^{\frac{3}{2}}}\right),
\end{array}
\end{displaymath}
Then, realizing similar operations for the other three couples
and putting together all the results we have, up to a $N^{-3/2}$ term,
\begin{multline*}
M_1=\frac{(p-1)!}{N^{p-1}} (p-1) 
\beta^2 u_N^2 
q^{p-2} (1-2q + \hat q_4)\\
\times
\sum_{{\J}'\in \asp} \sum_{{\J}\in \asp}  
\nu_{p-1,0}\Big((\eta^1_{\J}-\eta_{\J}^5)(\eta_{\J}^3-\eta^6_{\J})
(\eta_{{\J}'}^1-\eta_{{\J}'}^2)(\eta_{{\J}'}^{3}-\eta_{{\J}'}^{4})\Big)
.
\end{multline*}
Finally, from Corollary \ref{noudos}, (\ref{e6.s12}), 
and Corollary \ref{c3.6}, we can obtain
\begin{equation}\label{e6.21}
M_1=\beta^2 \frac{p(p-1)}{2}
q^{p-2} (1-2q +\hat q_4)
\nu\left(S_{1,2}^2\right)+
\mathcal{O}\left(\frac{1}{N^{\frac{3}{2}}}\right).
\end{equation}
Then, putting together (\ref{e6.11})-(\ref{e6.21}) we finish the proof.
\hfill
$\Box$

%\vspace{0.2cm}

It will be useful in the sequel to have some information about the 
correlations between $T_{l,l^{'}}$, $T_l$ and $T$. This is easily obtained in the following Proposition:
\begin{proposition}\label{p6.2}
The following cancellations hold true.
\begin{enumerate}
\item
If $l<l'$ and $(l,l')\neq (1,2)$, we have
$$
\nu\left(T_{1,2}\ T_{l,l'}\right)\ =\ 0.
$$
\item
For any $l$, we have
$$
\nu\left(T_{1,2}\ T_{l}\right)= 0, \qquad
\nu\left(T_{1,2}\ T\right) = 0.
$$
\item
For any $l \neq 1$, we have
$$
\nu\left(T_{1}\ T_{l}\right) = 0, \qquad
\nu\left(T_{1}\ T\right) = 0.
$$
\end{enumerate}
\end{proposition}

%\vspace{0.2cm}

\noindent
{\bf Proof:}
This is trivially obtained by some symmetry considerations among the different copies of $G_N$.
\hfill
$\Box$

Let us turn now to the second moment estimate of $T_1$.
\begin{proposition}\label{p6.3}
Whenever $\beta$ satisfies {\rm (H)}, we have
\begin{equation}\label{e6.301}
\Big|\nu\big(T_{1}^2\big)-\frac{B^2}{N}\Big|\le \frac{K}{N^{\frac{3}{2}}}
\end{equation}
where
\begin{equation} \label{e6.302}
B^2= (p-1)q^{p-2}(q-\hat q_4)\left[
\frac{(p-1)q^{p-2} + \beta^2\frac{p}{2}A^2}
{1-\beta^2 \frac{p(p-1)}{2} q^{p-2}(1-4q+3\hat q_4)}\right]. 
\end{equation}
\end{proposition}

%\vspace{0.2cm}

\noindent {\bf Proof:} 
We can decompose  $T_1$ as follows
$$
T_{1}=S_{1}+V_{N,T}^{1},
$$
with
\begin{eqnarray*}
S_{1}&=&
\frac{(p-1)!}{N^{p-1}}
\sum_{{\J} \in A_N^{p-1}}
(\eta^1_{\J}-b_{\J}) b_{\J},\label{e6.304}\\
V_{N,T}^{1}&=&
\frac{1}{N^{p-1}}
\sum_{{\J} \in \bar N_{p-1}^c}
(\eta^1_{\J}-b_{\J}) b_{\J}.\label{e6.305}
\end{eqnarray*}
The same kind of arguments as in Proposition \ref{p6.1} allow to state that
\begin{eqnarray}
\nu\left(S_{1}^2\right)&=&
\mathcal{O}\left(\frac{1}{N}\right),\label{e6.306} \\
\nu\left(T_{1}^2\right)&=&\nu\left(S_{1}^2\right)+
\mathcal{O}\left(\frac{1}{N^{\frac{3}{2}}}\right),\label{e6.307}
\end{eqnarray}
and
\begin{equation}\label{e6.308}
\nu\left(S_{1}^2\right)\ =\ \tilde M_1\ +\ \tilde M_2\ + \ \tilde M_3 +\mathcal{O}\left(\frac{1}{N^{\frac{3}{2}}}\right),
\end{equation} 
with
\begin{eqnarray*}
\tilde M_1&=&\frac{(p-1)!}{N^{p-1}} 
\nu\Bigg[(\ep^1-\ep^2) \ep^3
\sum_{{\J} \in \asp}
(\eta^1_{\J}-\eta_{\J}^4) \eta_{\J}^5\Bigg],\\
\tilde M_2&=&\frac{(p-1)!}{N^{p-1}}
\nu\Bigg[(\ep^1-\ep^2) \ep^3 
\sum_{{\J} \in \tilde Q_{N,p-1}^{p-1}}
(\eta^1_{\J}-\eta_{\J}^4) \eta_{\J}^5\Bigg],\\
\tilde M_3&=&\frac{(p-1)!}{N^{p-1}}
\nu\Bigg[(\ep^1-\ep^2) \ep^3 
\sum_{{\J} \in \bar Q_{N,p-1}^{p-1}}
(\eta^1_{\J}-\eta_{\J}^4) \eta_{\J}^5\Bigg].\\
\end{eqnarray*}
By Lemma \ref{cardbn} we  easily get
\begin{equation}\label{e6.309}
\big|\tilde M_3\big| \le \frac{L}{N^2}.
\end{equation}
Using exactly the same arguments as in the study of 
$M_2$, we obtain 
\begin{equation}\label{e6.310}
\tilde M_2= \frac{1}{N}(p-1)^2 q^{2(p-2)} (q - \hat q_4)     
+\mathcal{O}\left(\frac{1}{N^{\frac{3}{2}}}\right).
\end{equation}
Finally, we  deal with $\tilde M_1$. Let
\begin{displaymath}
f=
\frac{(p-1)!}{N^{p-1}}(\ep^1-\ep^2) \ep^3
\sum_{{\J} \in \asp}
(\eta^1_{\J}-\eta_{\J}^4) \eta_{\J}^5.
\end{displaymath}
Then, Proposition \ref{p5.1}, Corollary \ref{c3.6} and one symmetry consideration yield
\begin{eqnarray*}
\tilde M_1&=&  \  
\Big[\nu_{p-1,0}(f)+ \nu'_{p-1,0}(f)\Big] 
+\mathcal{O}\left(\frac{1}{N^{\frac{3}{2}}}\right)\\
&=&
 \nu'_{p-1,0}(f)  \ 
+\mathcal{O}\left(\frac{1}{N^{\frac{3}{2}}}\right).
\end{eqnarray*}
Using Proposition \ref{pderivnu} we can write
$$
\nu_{p-1,0}'(f)
=\beta^2u_N^2\sum_{J'\in Q_{N,p-1}^p}
\big(\tilde F^{J'}_1+\tilde F^{J'}_2+ \tilde F^{J'}_3\big),
$$
with
\begin{eqnarray*}
\tilde F_1^{J'}&=& 
\nu_{p-1,0}\Bigg( f
\sum_{1\le l<l'\le 5}(\eta_{J'}^l\eta_{J'}^{l'}-q^{p-1}) 
\ep_{J'}^l\ep_{J'}^{l'}\Bigg),\\
\tilde F_2^{J'}&=&-5 \nu_{p-1,0}\Bigg( f
\sum_{l\le 5}(\eta_{J'}^l\eta_{J'}^{6}-q^{p-1}) 
\ep_{J'}^l\ep_{J'}^{6}\Bigg),\\
\tilde F_3^{J'}&=&15 \nu_{p-1,0}\Big( f
(\eta_{J'}^6\eta_{J'}^{7}-q^{p-1}) 
\ep_{J'}^6\ep_{J'}^{7}\Big).
\end{eqnarray*}
The only non-vanishing terms in these expressions are the one induced by the couples in $W_1$, where
$$W_1=\{(1,3), (2,3), (1,4), (2,4), (1,5), (2,5), (1,6), (2,6)\},$$ 
We can then rewrite $\tilde M_1$ as follows:
$$\tilde M_1= \sum_{(l,l')\in W_1}
\tilde F_{(l,l')}
+\mathcal{O}\left(\frac{1}{N^{\frac{3}{2}}}\right),$$
with
\begin{eqnarray*}
\tilde F_{(l,l')} 
&=&\beta^2 u_N^2 
\sum_{J'\in Q_{N,p-1}^p}
\nu_{p-1,0}\Bigg( f
(\eta_{J'}^l\eta_{J'}^{l'}-q^{p-1}) 
\ep_{J'}^l\ep_{J'}^{l'}\Bigg),\quad  l'\neq 6,\\
\tilde F_{(l,6)}&=&-5 
\beta^2 u_N^2  \sum_{J'\in Q_{N,p-1}^p}   
\nu_{p-1,0}\Bigg( f
(\eta_{J'}^l\eta_{J'}^{6}-q^{p-1}) 
\ep_{J'}^l\ep_{J'}^{6}\Bigg).
\end{eqnarray*}
We first analize together the couples $(1,6), (2,6)$. The most important
remark is the following consequence of Lemma \ref{l2.5}:
$$ 
\nu_{p-1,0}\Big((\ep^1-\ep^2) \ep^3 \ep_{p-1}^1\ep_{p-1}^{6}\Big)
=
-\nu_{p-1,0}\Big((\ep^1-\ep^2) \ep^3 \ep_{p-1}^2\ep_{p-1}^{6}\Big)
=q^{p-2}(q-\hat q_4).$$
Then, by means of the arguments used in Proposition \ref{p6.1} we can prove that
$$\tilde F_{(1,6)}  +\tilde F_{(2,6)} =
-5 \beta^2 \frac{p(p-1)}{2} q^{p-2}(q-\hat q_4) \nu\left(S_1^2\right)
 +\mathcal{O}\left(\frac{1}{N^{\frac{3}{2}}}\right).$$
Performing the same sort of computations, we can get
$$
\tilde F_{(1,3)}  +\tilde F_{(2,3)} =
 \beta^2 \frac{p(p-1)}{2} q^{p-2}(1-q) \nu\left(S_1^2\right)
 +\mathcal{O}\left(\frac{1}{N^{\frac{3}{2}}}\right).$$
We finish the study of $\tilde M_1$ by considering the last four couples. 
As usual, independence, symmetry, relations between the different
sets, Proposition \ref{p3.2}, Lemma \ref{l2.5} and  Corollary \ref{c3.6}
imply that, up to a term of order $N^{-3/2}$,
$$\begin{array}{l}
\displaystyle\tilde F_{(1,4)}  +\tilde F_{(2,4)} +
\tilde F_{(1,5)}  +\tilde F_{(2,5)} 
\displaystyle =
\beta^2 \frac{p(p-1)}{2} q^{p-2}(q-\hat q_4) 
     \left(\frac{(p-1)!}{N^{p-1}}\right)^2 \\
     \displaystyle \ \times
     \sum_{{\J}'\in \asp} \sum_{{\J}\in \asp}  
\nu\Big((\eta^1_{\J}-\eta_{\J}^4) \eta^5_{\J} 
(\eta_{{\J}'}^1-\eta_{{\J}'}^{2})(\eta_{{\J}'}^4+\eta_{{\J}'}^{5})\Big)
\end{array}$$
Then, (\ref{equasif}), (\ref{e6.306}), (\ref{e6.4}), 
Propositions \ref{p6.2}, 
 \ref{p3.2}, \ref{p5.1}, (\ref{e6.5}), (\ref{e6.6}) and (\ref{e6.307}) yield
$$\begin{array}{l}
\displaystyle\tilde F_{(1,4)}  +\tilde F_{(2,4)} +
\tilde F_{(1,5)}  +\tilde F_{(2,5)}\\ 
\quad \displaystyle
=\beta^2 \frac{p(p-1)}{2} q^{p-2}(q-\hat q_4) 
\nu\Bigg(\Big(R_{1,5}^{p-1}-    
R_{4,5}^{p-1}\Big)\\
\quad\qquad \displaystyle \quad\times
\Big(R_{1,4}^{p-1} -R_{2,4}^{p-1}-R_{2,5}^{p-1}+  
R_{1,5}^{p-1}\Big)\Bigg)  
 +\mathcal{O}\left(\frac{1}{N^{\frac{3}{2}}}\right)\\
\quad \displaystyle
=\beta^2 \frac{p(p-1)}{2} q^{p-2}(q-\hat q_4) 
\nu\Bigg(\Big(T_{1,5}-T_{4,5}+T_1-T_4\Big)\\
\quad \displaystyle\qquad \quad\times
\Big(T_{1,4}-T_{2,4}-T_{2,5}+T_{1,5}
+2T_1-2T_2\Big)\Bigg)
 +\mathcal{O}\left(\frac{1}{N^{\frac{3}{2}}}\right)\\
\quad \displaystyle
=\beta^2 \frac{p(p-1)}{2} q^{p-2}
(q-\hat q_4)  \Bigg(\frac{A^2}{N}+ 2\nu\left(S_1^2\right)\Bigg)
 +\mathcal{O}\left(\frac{1}{N^{\frac{3}{2}}}\right).
\end{array}$$
So, 
\begin{equation}\label{e6.311}
\tilde M_1=
 \beta^2 \frac{p(p-1)}{2} q^{p-2}
 \Big[(1-4q+3\hat q_4) \nu\left(S_1^2\right)
+ (q-\hat q_4)  \frac{A^2}{N}\Big]
 +\mathcal{O}\left(\frac{1}{N^{\frac{3}{2}}}\right).
 \end{equation}
Finally, putting together (\ref{e6.307})-(\ref{e6.311}) we get (\ref{e6.301}) and (\ref{e6.302}).

\hfill
$\Box$

%\vspace{0.2cm}

The second moment of $T$ can also be estimated in the following way:
\begin{proposition}\label{p6.5}
If $\beta$ satisfies {\rm (H)}, we have
$$
\Big|\nu\big(T^2\big)-\frac{C^2}{N}\Big|\le \frac{L}{N^{\frac{3}{2}}}
$$
where
\begin{multline*}
C^2=\\
(p-1)q^{p-2}
\left[\frac{(\hat q_4-q^2)\big[(p-1) q^{p-2}+\beta^2\frac{p}{2} A^2\big] 
+\beta^2 p (2q+q^2-3\hat q_4)B^2}
{1-\beta^2 \frac{p(p-1)}{2} q^{p-2}(1-4q+3\hat q_4)}\right].
\end{multline*}
\end{proposition}

%\vspace{0.2cm}

\noindent {\bf Proof:}
The proof of this relation goes along the same lines as Propositions  \ref{p6.1} and \ref{p6.3}, and we will only point out the main differences with the latter.

Observe that
\begin{equation}\label{e6.508}
\nu\left(T^2\right)\ =\ \hat M_1\ +\ \hat M_2\ 
+\mathcal{O}\left(\frac{1}{N^{\frac{3}{2}}}\right),
\end{equation} 
with
\begin{eqnarray*}
\hat M_1&=&\frac{(p-1)!}{N^{p-1}} 
\nu\Bigg[(\ep^1 \ep^2 -q^{p-1}) 
\sum_{{\J}\in \asp}
(\eta^3_{\J}\eta_{\J}^4-q^{p-1})\Bigg],\\
\hat M_2&=&\frac{(p-1)!}{N^{p-1}}
\nu\Bigg[(\ep^1 \ep^2 -q^{p-1}) 
\sum_{{\J} \in \tilde Q_{N,p-1}^{p-1}}
(\eta^3_{\J} \eta_{\J}^4-q^{p-1})\Bigg].
\end{eqnarray*}
The same computation as in Proposition \ref{p6.3} yield, for the term $\tilde M_2$, the following equality: 
\begin{equation}\label{e6.510}
\hat M_2= \frac{1}{N}(p-1)^2 q^{2(p-2)} (\hat q_4 - q^2)     
+\mathcal{O}\left(\frac{1}{N^{\frac{3}{2}}}\right).
\end{equation}
In order to deal with $\hat M_1$ 
we have to evaluate the derivative 
of $\nu_{p-1,0}(f)$,
where
\begin{displaymath}
f=
\frac{(p-1)!}{N^{p-1}}(\ep^1\ep^2-q^{p-1}) 
\sum_{{\J} \in \asp}
(\eta^3_J \eta_{\J}^4-q^{p-1}).
\end{displaymath}
As in Proposition \ref{p6.3}, by arguments of symmetry, we can rewrite
$\hat M_1$ as
\begin{equation}\label{e6.511}
\hat M_1= 
\sum_{(l,l')\in W_2} \hat F_{(l,l')}
+\mathcal{O}\left(\frac{1}{N^{\frac{3}{2}}}\right),
\end{equation}
where
$$W_2=\{(1,2), (1,3), (2,3), (1,4), (2,4), (3,4),
(1,5), (2,5), (3,5), (4,5), (5,6)\},$$ 
and, for $(l,l')\in W_2$,

\begin{eqnarray*}
\hat F_{(l,l')}&=& \beta^2 u_N^2 (p-1)!
\sum_{J'\in Q_{N,p-1}}
\nu_{p-1,0}\Bigg( f
(\eta_{J'}^l\eta_{J'}^{l'}-q^{p-1}) 
\ep_{J'}^l\ep_{J'}^{l'}\Bigg),\,  l'\neq 5,6,\\
\hat F_{(l,5)}&=&-4 
\beta^2 u_N^2 (p-1)! \sum_{J'\in Q_{N,p-1}}   
\nu_{p-1,0}\Bigg( f
(\eta_{J'}^l\eta_{J'}^{5}-q^{p-1}) 
\ep_{J'}^l\ep_{J'}^{5}\Bigg),\\
\hat F_{(5,6)}&=&10
\beta^2 u_N^2 (p-1)! \sum_{J'\in Q_{N,p-1}}   
\nu_{p-1,0}\Bigg( f
(\eta_{J'}^5\eta_{J'}^{6}-q^{p-1}) 
\ep_{J'}^5\ep_{J'}^{6}\Bigg).\\
\end{eqnarray*}
Let
$$K_{\beta,p,q}=\beta^2 \frac{p(p-1)}{2}q^{p-2}. 
$$  
Operating as in Propositions \ref{p6.1} and \ref{p6.3}, we can obtain
\begin{equation}\label{e6.512}
\begin{array}{l}
\displaystyle\hat F_{(1,2)}= K_{\beta,p,q}\ (1-q^2) 
\nu\left(T^2\right),\\[3mm]
\displaystyle\hat F_{(l,l')}= K_{\beta,p,q}\ q (1-q) 
[\nu\left(T^2\right)+B^2],\ {\rm for } (l,l')=
(1,3), (1,4), (2,3), (2,4),\\[3mm]
\displaystyle\hat F_{(3,4)}= K_{\beta,p,q}\ (\hat q_4-q^2) 
[\nu\left(T^2\right)+2B^2+A^2],\\[3mm]
\displaystyle\hat F_{(l,l')}=-4 K_{\beta,p,q}\ q (1-q) 
\nu\left(T^2\right),\ {\rm for } (l,l')=
(1,5),  (2,5),\\[3mm]
\displaystyle\hat F_{(l,l')}=-4 K_{\beta,p,q}\  (\hat q_4-q^2) 
[\nu\left(T^2\right)+B^2],\quad{\rm for }\ (l,l')=
(3,5), (4,5),\\[3mm]
\displaystyle\hat F_{(5,6)}=10  K_{\beta,p,q}\  (\hat q_4-q^2) 
\nu\left(T^2\right).
\end{array}\end{equation}
Then, from  (\ref{e6.508})-(\ref{e6.512})  
we can conclude the proof of this proposition.

\hfill
$\Box$

%\vspace{0.1cm}

\section{Central Limit Theorems}

The main result of this section will be a CLT for the fluctuations of $R_{1,2}$, though we will get, on our way to the proof of this theorem, some general limit relations for the joint fluctuations of $T_{l,l^{'}}$, $T_l$ and $T$. First, observe that, as an immediate consequence of Propositions \ref{p6.1}, 
\ref{p6.2}, \ref{p6.3}, \ref{p6.5} and equality (\ref{e6.4})
we have the following
\begin{proposition}\label{p7.1}
If $\beta$ verifies {\rm (H)}, we have
$$\left|\nu \Big( (R_{1,2}^{p-1}-q^{p-1})^2 \Big)
-\frac{1}{N}\left(A^2+2B^2+C^2\right)\right| 
\le \frac{K}{N^{\frac{3}{2}}}.$$
\end{proposition}
Our aim here will be to generalize this estimate, 
obtaining a similar relation for 
$$\nu\left(\left(R_{1,2}^{p-1}-q^{p-1} \right)^k\right).$$

Let us first state the result we obtain for the fluctuations of $T_{l,l^{'}}$: define $a(k)=\be g^k$, for 
a standard Gaussian random variable $g$.

\begin{theorem}\label{t7.2}
Let $\beta$ satisfying {\rm (H)}, $n\in \mathbb{N}$. For any couple $1\le l<
l'\le n$,  consider an integer $k(l,l')\ge 0$.
Set $k=\sum_{l<l'} k(l,l') $.
Then, we have
\begin{displaymath}
\left|\nu\left(\prod_{l<l'} T_{l,l'}^{k(l,l')} \right)-
\frac{1}{N^\frac{k}{2}}\prod_{l<l'} a(k(l,l'))  A^k\right|
\le \frac{L(k)}{N^\frac{k+1}{2}}.
\end{displaymath}
\end{theorem}

%\vspace{0.1cm}

\noindent {\bf Proof:} 
We will follow the proof of Theorem 2.7.1 in \cite{Tbk}, and
we prove this theorem by induction over $k$. 
The case $k=2$ has been proved in Section 6, Proposition \ref{p6.1}.

We assume now that the result is true up to the order $k-1$. Then, using  (\ref{e6.7}) 
and (\ref{YY}), we can easily check that if $\sum_{l<l'} k(l,l')=j$ where $j \le k-1$,
we have
\begin{equation}\label{indu}
\nu\left(\prod_{l<l'} S_{l,l'}^{k(l,l')} \right)=
\mathcal{O}\left(\frac{1}{N^\frac{j}{2}} \right).
\end{equation}

Furthermore, using the induction, (\ref{indu}), and relations (\ref{e6.7}) and (\ref{YY}), we also have
\begin{equation}\label{e7.0}
\nu\left(\prod_{l<l'} T_{l,l'}^{k(l,l')}\right) =
\nu\left(\prod_{l<l'} S_{l,l'}^{k(l,l')}\right)
+\mathcal{O}\left(\frac{1}{N^\frac{k+1}{2}} \right).
\end{equation}

Observe that $\prod_{l<l'} S_{l,l'}^{k(l,l')}$ can be decomposed as 
\begin{displaymath}
\prod_{l<l'} S_{l,l'}^{k(l,l')} =
\prod_{1\le v \le k} S_{l(v),l'(v)},
\end{displaymath}
where, for any integer $v \le k$, $l(v), l'(v)$ are two integers such that
$$
(l(v),l'(v))=(1,2)
\Longleftrightarrow v\le k(1,2).
$$
We can assume without loss of generality that   $k(1,2)\ge 1$.
We consider, for $1\le v \le k$, integers $j(v), j'(v)$, all
different, and greater than $n$. Thus
$$
\nu\left(\prod_{l<l'} S_{l,l'}^{k(l,l')}\right) =
\nu\left(\prod_{1\le v \le k} S_{l(v),l'(v)}\right)
=  \nu\left(\prod_{1\le v \le k} U(v)\right),
$$
with
\begin{displaymath}
U(v) =
\frac{(p-1)!}{N^{p-1}} 
\left(\hat T\left(\eta^{l(v)}\right)-
\hat T\left(\eta^{j(v)}\right)\right)\cdot   
\left(\hat T\left(\eta^{l'(v)}\right)-   
\hat T\left(\eta^{j'(v)}\right)\right).   
\end{displaymath}
Set
$$\ep(v)=\Big(\ep^{l(v)}-\ep^{j(v)}\Big)
\left(\ep^{l'(v)}-\ep^{j'(v)}\right).$$
By the usual symmetry argument, and using (\ref{indu}) and Lemma \ref{cardasp}  we have
\begin{equation}\label{e7.3}
\nu\left(\prod_{l<l'} S_{l,l'}^{k(l,l')}\right) =
\nu\left(\ep(1) \prod_{2\le v \le k} U(v)\right)
+\mathcal{O}\left(\frac{1}{N^\frac{k+1}{2}} \right).
\end{equation}
For any $v\ge k$, $U(v)$ can be written as
$$U(v)=U_1(v)+U_2(v)+U_3(v),$$
where $U_1(v), U_2(v), U_3(v)$ are defined by means of 
$\asp,  \tilde Q_{N,p-1}^{p-1},   \bar Q_{N,p-1}^{p-1}$,
respectively.
Similarly to Proposition \ref{p6.1}, we can get
\begin{eqnarray}\label{e7.4}
\nu\left(\ep(1) \prod_{2\le v \le k} U(v)\right)&=&
\nu\left(\ep(1) \prod_{2\le v \le k} \left(U_1(v)+U_2(v)\right)\right)
+\mathcal{O}\left(\frac{1}{N^\frac{k+1}{2}} \right)\nonumber\\
&=& \nu\left(\ep(1) \prod_{2\le v \le k} U_1(v)\right)
+ I 
+\mathcal{O}\left(\frac{1}{N^\frac{k+1}{2}} \right),
\end{eqnarray}
with
\begin{displaymath}
I=\sum_{2\le u\le k} \nu\left(\ep(1) U_2(u)
\prod_{v\neq u} U_1(v)\right),
\end{displaymath}
where $\prod_{v\neq u}$ means that the product is over 
$2\le v\le k$, $v\neq u$.

We now study $I$. As in \cite{Tbk},
Corollary \ref{c3.6} and the  usual procedure
imply, if $k(1,2)\ge2$,
\begin{equation}\label{e7.5}
I=\frac{(k(1,2)-1)}{N}(p-1)^2
q^{2(p-2)}(1-2q+\hat q_4)
\nu\left(\prod_{3\le v \le k} U_1(v)\right)
+\mathcal{O}\left(\frac{1}{N^\frac{k+1}{2}} \right).
\end{equation}
If $k(1,2)=1$, $I=0$.

Now we deal with the other term of (\ref{e7.4}). Using again
Corollary \ref{c3.6}, we should study
particularly the derivative of
$$\nu_{p-1,0}\left(\ep(1) \prod_{2\le v \le k} U_1(v)\right).$$ 
Indeed, the following relation is not difficult to obtain:
\begin{equation}\label{e7.6}
\begin{array}{l}
\displaystyle\nu \left(\ep(1) \prod_{2\le v \le k} U_1(v)\right)\\[3mm]
\qquad = \displaystyle\nu_{p-1,0}'\left(\ep(1) \prod_{2\le v \le k} U_1(v)\right)
+\mathcal{O}\left(\frac{1}{N^\frac{k+1}{2}} \right)\\[3mm]
\qquad\displaystyle=
\beta^2\frac{p(p-1)}{2}
q^{p-2}(1-2q+\hat q_4)
\nu\left(U_1(1)\prod_{2\le v \le k} U(v)\right)
+\mathcal{O}\left(\frac{1}{N^\frac{k+1}{2}} \right).
\end{array}\end{equation}
Putting together (\ref{e7.0}), (\ref{e7.3})-(\ref{e7.6})
and  reasoning by induction over $k$ 
we can conclude the proof of this proposition.
\hfill
$\Box$

%\vspace{0.1cm}

In a similar way to the previous theorem we can prove the following
relations on the joint behavior  of $T_{l,l'}$, $T_l$ and $T$ (we do not include the proofs here, since they follow closely the lines of \cite{Tbk}).

\begin{theorem}\label{t7.3}
Let $\beta$ satisfying {\rm (H)}, $n\in \mathbb{N}$. 
For $1\le l<  l'\le n$,  consider integers $k(l,l')\ge 0$ and 
$k_1=\sum_{1 \le l<l'\le n} k(l,l')$. 
For $1\le l\le n$,  let $k(l)$ be a positive integer and set 
$k_2=\sum_{1 \le l \le n} k(l)$. 
Then, if $k=k_1+k_2$, we have
\begin{multline*}
\Bigg|\nu\left(\prod_{1 \le l<l'\le n} T_{l,l'}^{k(l,l')} 
\prod_{1 \le l\le n} T_{l}^{k(l)}
\right)\\
-
\frac{1}{N^\frac{k}{2}}\prod_{1 \le l<l'\le n} a(k(l,l'))  
\prod_{1 \le l\le n} a(k(l)) A^{k_1} B^{k_2}\Bigg|
\le \frac{L(k)}{N^\frac{k+1}{2}}.
\end{multline*}
\end{theorem}

%\vspace{0.1cm}

\begin{theorem}\label{t7.4}
Let $\beta$ satisfying {\rm (H)}, $n\in \mathbb{N}$. 
For $1\le l<  l'\le n$,  consider integers $k(l,l')\ge 0$ and 
$k_1=\sum_{1 \le l<l'\le n} k(l,l')$. 
For $1\le l\le n$,  let $k(l)\ge 0$ and 
$k_2=\sum_{1 \le l\le n} k(l)$. Let $k_3\in \mathbb{N}$.  
Then, if $k=k_1+k_2+k_3$, we have
\begin{displaymath}\begin{array}{l}
\displaystyle\Bigg|\nu\left(\prod_{1 \le l<l'\le n} T_{l,l'}^{k(l,l')} 
\prod_{1 \le l\le n} T_{l}^{k(l)} T^{k_3}
\right)\\[4mm]
\displaystyle\quad\quad-
\frac{1}{N^\frac{k}{2}}\prod_{1  \le l<l'\le n} a(k(l,l'))  
\prod_{1 \le l\le n} a(k(l)) a(k_3) A^{k_1} B^{k_2} C^{k_3}\Bigg|
\le \frac{L(k)}{N^\frac{k+1}{2}}.
\end{array}\end{displaymath}
\end{theorem}

%\vspace{0.1cm}

All the preceding considerations allow us to get the following  Central Limit Theorem.

\begin{theorem}\label{t7.5}
Let $\beta$ satisfying {\rm (H)}, $\hat k\in \mathbb{N}$. 
Then,
$$
\left|\nu\left(\left(R_{1,2}^{p-1}-q^{p-1}\right)^{\hat k}\right)
-\frac{1}{N^\frac{\hat k}{2}}\ a(\hat k) \left(A^2+2B^2+C^2\right)^\frac{\hat k}{2}
\right|
\le \frac{L(\hat k)}{N^{\frac{\hat k+1}{2}}}.
$$
\end{theorem}

%\vspace{0.1cm}

\noindent {\bf Proof:} 
It is well-known that for any $\tilde k \in \mathbb{N}$,
\begin{equation}\label{e7.8}\begin{array}{l}
\displaystyle a(2\tilde k)= 
\displaystyle\frac{(2\tilde k)!}{2^{\tilde k}\tilde k!},\\[3mm]
\displaystyle a(2\tilde k+1)=0.  
\end{array}
\end{equation}

By (\ref{e6.4}),   a combinatorial property,
Theorems \ref{t7.2}, \ref{t7.3}
and \ref{t7.4} and  (\ref{e7.8}), we have
\begin{eqnarray*}
&&\nu\left(\left(R_{1,2}^{p-1}-q^{p-1}\right)^{\hat k}\right)
= \nu\left(\left(T_{1,2}+T_1+T_2+T\right)^{\hat k}\right)\\ 
&&=\sum
\frac{\hat k!}{k_{1,2}! k_1! k_2! k!}
\nu\left(T_{1,2}^{k_{1,2}}\ T_1^{k_1}\ T_2^{k_2} \
T^{k}\right)\\
&&=\frac{1}{N^\frac{\hat k}{2}}\sum
\frac{a(k_{1,2})\ a(k_1)\ a(k_2)\ a(k)\ \hat k!}{k_{1,2}! k_1! k_2! k!}
A^{k_{1,2}}\ B^{k_1+k_2}\ C^{k} 
+\mathcal{O}\left(\frac{1}{N^\frac{\hat k+1}{2}}\right)\\
&&=\frac{1}{N^\frac{\hat k}{2}}\sum_{\rm even}
\frac{a(\hat k) \left(\frac{\hat k}{2}\right)!
\left(A^2\right)^\frac{k_{1,2}}{2} 
 \left(B^2\right)^\frac{k_1+k_2}{2} 
 \left(C^2\right)^\frac{k}{2}}
{\left(\frac{k_{1,2}}{2}\right)! 
\left(\frac{k_1}{2}\right)!
\left(\frac{k_2}{2}\right)!
\left(\frac{k}{2}\right)! } 
+\mathcal{O}\left(\frac{1}{N^\frac{\hat k+1}{2}}\right)\\
&&=\frac{1}{N^\frac{\hat k}{2}}a(\hat k) \left(A^2+2B^2+C^2\right)^\frac{\hat k}{2}
+\mathcal{O}\left(\frac{1}{N^\frac{\hat k+1}{2}}\right),
\end{eqnarray*}
where $\sum$ means the summatory of 
$k_{1,2}, k_1, k_2, k\in 
\mathbb{N}$ such that 
$k_{1,2}+k_1+k_2+k=\hat k$ and $\sum_{\rm even}$ means the summatory of 
$k_{1,2}, k_1, k_2, k\in 
\mathbb{N}$ such that all these numbers are even and
$k_{1,2}+k_1+k_2+k=\hat k$. 
\hfill
$\Box$

%\vspace{0.1cm}

Eventually, a CLT for $R_{1,2}$, that can be considered as the main result of this section, is easily obtained from the last theorem.
\begin{corollary}\label{c7.9}
Let $\beta$ satisfying {\rm (H)}, $\hat k\in \mathbb{N}$. 
Then, we have
\begin{equation}\label{et7.9}
\left|\nu\left(\left(R_{1,2}-q\right)^{\hat k}\right)
-\frac{1}{N^\frac{\hat k}{2}}\ a(\hat k) \frac{\left(A^2+2B^2+C^2\right)^\frac{\hat k}{2}}
{((p-1)q^{p-2})^{\hat k}}
\right|
\le \frac{L(\hat k)}{N^{\frac{\hat k+1}{2}}}.
\end{equation}
\end{corollary}

%\vspace{0.2cm}

\noindent {\bf Proof:} 
>From (\ref{emajoration})
and since $q$ is a strictly positive number we have
\begin{displaymath}
(R_{1,2}-q)=\frac{1}{(p-1) q^{p-2}}(R_{1,2}^{p-1}-q^{p-1})
-\frac{p-2}{2q^{p-2}}\ \xi^{p-3} (R_{1,2}-q)^2,
\end{displaymath}
where $\xi \in  (R_{1,2} \wedge q, R_{1,2} \vee q)$.

Since $\xi\le 1$ 
we obtain (\ref{et7.9})
by means of  Theorems \ref{t7.5} and \ref{t3.3}.    
\hfill
$\Box$

\section{Appendix}

In this appendix we will recall the definitions of all the sets appearing
in the paper, as well as some results about their size that will be used
throughout this paper. Since the method and the tools needed to prove these results are always the same, we will only give some examples.
Recall that $P_{m}(N)$ denotes a  polynomial of order $m$ in $N$.

\begin{definition}\label{setajr}
For $w\ge 1$ and $1 \le r \le w$, set
$$
A_w^r = \lcl (i_1,\ldots,i_r)\in\mathbb{N}^r; 1\le i_1<\cdots<i_r\le w  
\rcl.
$$
\end{definition}

\begin{lemma}\label{cardasp}
For $N\ge 1$ and $p\ge 1$, we have
$$
|A_{N}^p| = {N \choose p} = \frac{N^p}{p!} + P_{p-1}(N).$$
\end{lemma}

\begin{definition}\label{setqmjr}
For $w\ge 1$ and $1 \le r \le w$, set
$$
Q_{w,j}^r= \big\{(i_1,\dots,i_r)\in \mathbb{N}^r;1 \le  i_1<\cdots <i_r \le w, i_r>w-j\big\}.
$$
The set $Q_{w,j}^r$ can be split into
$$Q_{w,j}^r= \bar Q_{w,j}^r  \cup \tilde Q_{w,j}^r,$$
where
\begin{eqnarray*}
\bar Q_{w,j}^r&=&\big\{(i_1,\dots,i_r)\in Q_{w,j}^r\ ;\ i_1<\cdots<i_r, i_{r-1}>w-j\big\},\\
\tilde Q_{w,j}^r&=&\big\{(i_1,\dots,i_r)\in Q_{w,j}^r\ ;\ i_1<\cdots<i_r, i_{r-1} \le w-j \big\}.
\end{eqnarray*}
\end{definition}

\begin{lemma}\label{cardbn}
For $N\ge 1$ and $k\ge 1$, we have
\begin{eqnarray}
(p-1)! |Q_{N,1}^p|&=& N^{p-1}-\frac{p(p-1)}{2}N^{p-2}+P_{p-3}(N),\nonumber\\
(p-1)! |Q_{N,k}^p|&=& k N^{p-1}+P_{p-2}(N),\nonumber\\
|\bar Q_{N,k}^p|&=& P_{p-2}(N),\nonumber\\
\big| \tilde Q_{N,p}^{p}\big| &=&P_{p-1}(N).\nonumber
\end{eqnarray}
\end{lemma}

%\vspace{0.1cm}

\noindent
{\bf Proof:} 
It is easily seen that
\begin{eqnarray*}
(p-1)! |Q_{N,1}^p|&=&
(N-1)\cdots(N-(p-1))\\
&=&\displaystyle N^{p-1}-\lp\sum_{j=1}^{p-1}j\rp N^{p-2}
+P_{p-3}(N),
\end{eqnarray*}
which implies our first claim. In order to prove the second one, we use
the following fact
\begin{eqnarray*}
|Q_{N,k}^p|&=&{N \choose p} - {N-k \choose p}\\
&=&\displaystyle\frac{1}{p!} \Big(\sum_{j=k}^{k+p-1}j 
- \sum_{j=1}^{p-1}j\Big) N^{p-1} + P_{p-2}(N)\\
&=& \frac{k}{(p-1)!} N^{p-1} + P_{p-2}(N).
\end{eqnarray*}
Finally
$$|\bar Q_{N,k}^p|={N \choose p}-{N-k \choose p}-k{N-k \choose p-1}
=P_{p-2}(N),$$
and this finishes the proof of the third claim. The last one 
is an easy consequence of the previous results.
\hfill $\Box$

\begin{definition}\label{setnr}
For $r\ge 1$, set
$$
N_{r}=
\{(i_1,\dots,i_{r})\in \{1,\dots,N\}^{r}\}.
$$
The set $N_{r}$ can be split into
$$N_{r}= \bar N_{r}  \cup \bar N_{r}^c,$$
where
$$
\bar N_{r}= \big\{(i_1,\dots,i_r)\in N_{r}\ ;\ i_j \neq i_k {\rm \, for \,  all \, } j \neq k \big\}.
$$
\end{definition}

\begin{lemma}\label{cardnp}
For $N\ge 1$ and $p\ge 2$, we have
$$
|\bar N^c_{p-1}|=\frac{(p-1)(p-2)}{2} N^{p-2} + P_{p-3}(N).
$$
\end{lemma}

%\vspace{0.1cm}

\end{document}